\documentclass[11pt,leqno,twoside]{article} 
\usepackage{amsmath,amssymb,fancyhdr} 
\usepackage{bm} 

\newcommand{\ovl}[1]{\hskip1pt\raise10pt\hbox{\leaders\hrule \hskip7pt}\hskip-9pt#1}
\newcommand{\ovls}[1]{\hskip1pt\raise7.425pt\hbox{\leaders\hrule \hskip7pt}\hskip-9pt#1}
\newcommand{\wovl}[1]{\hskip1pt\raise10pt\hbox{\leaders\hrule \hskip8.5pt}\hskip-9.75pt#1}
\def\wtilde{\widetilde}
\def\what{\widehat}
\def\wbar{\overline}

\newcommand{\op}[1]{\operatorname{\text{\rm #1}}} 

\def\sgn{\op{sign}}\def\spn{\op{span}}\def\graph{\op{graph}}
\def\dist{\op{dist}} \def\dvg{\op{div}} 
\def\spt{\op{support}}

\def\R{\mathbb{R}} \def\C{\mathbb{C}}   
\def\chir{\raise2pt\hbox{$\chi$}}
\def\ha{{\textstyle\frac{1}{2}}}
\def\fr #1/#2{{\textstyle\frac{#1}{#2}}}  \def\sfr #1/#2{\text{\footnotesize$\frac{#1}{#2}$}} 
\def\Fr #1/#2{{\displaystyle\frac{#1}{#2}}}  

\newcommand{\prelemskip}{\vskip 5pt plus 1.5pt minus 1.5pt}
\newcommand{\postlemskip}{\vskip 4.5pt plus 1pt minus 0.5pt} 
\newenvironment{state}[1]{\par\prelemskip\noindent {\bf#1} %
\textit\bgroup\abovedisplayskip=5pt plus 3pt minus 1pt%
\belowdisplayskip=5pt plus 3pt minus 1pt}%
{\egroup \par\ifdim\lastskip<\medskipamount \removelastskip\penalty55\postlemskip\fi}  

\renewcommand{\epsilon}{\varepsilon}
\def\al#1{$$\begin{aligned}#1\end{aligned}$$}
\def\lal#1{\begin{align*}#1\end{align*}}
\def\mint{\text{\large$\int\hskip-3pt$}}  
  
\def\sint{\text{\small$\int\hskip-4pt$}}  
\def\tsum{{\textstyle\sum}}
\def\ssum{{\textstyle{\sum}}} 
\newcommand{\Bl}[1]{\text{\small$\Bigl#1$}} 
\newcommand{\Br}[1]{\text{\small$\Bigr#1$}}

\def\bn#1{\text{\bf #1}}

\DeclareMathSymbol{\pri}{\mathcal}{symbols}{"30}
\renewcommand{\prime}{\hskip0.8pt\pri\hskip-0.4pt}

\title{\vskip-.6in  A frequency function and singular set bounds\\ for branched minimal
immersions \vskip-.1in
 \author{\scshape Leon Simon\thanks{Partly supported by NSF DMS--0104049 \&
DMS-0406209 at Stanford University}\,\, \& Neshan Wickramasekera\thanks{\hskip0.3pt 
Partly supported by NSF DMS-0601265 \hskip0.5pt \& DMS-0707005 at U.C.~San Diego }
}\date{\vspace{-.55in}}}%

\begin{document}

\maketitle

\large
\renewcommand{\contentsname}{\centerline{\small\rm Contents}}


\bigskip

\centerline{Contents}

\medskip

{\par{\small   \multiply \baselineskip by 12 \divide \baselineskip by 13   
\halign{\hskip0.5in \S # &\hskip4pt \hbox to 4.3truein{  # \dotfill}  &\,\,\hfill # \cr 
{0} & Introduction & {1} \cr  
{1} &  Preliminaries & {3} \cr 
{2} &  2-valued $C^{1}$ harmonic functions---Part I & {6} \cr 
{3} &  $C^{1,\alpha }$ estimates for a class of linear equations & {12} \cr  
{4} &   2-valued $C^{1}$ harmonic functions---Part II &{15}   \cr 
{5} &  Regularity for $u_{a}= {\textstyle \frac {1}{2}} (u_{1}+u_{2})$ and %
                                                             $v=\pm \ha(u_{1}-u_{2})$, Part I &{19}  \cr    
{6} &  Some Growth Results for a Class of Linear Equations &{25} \cr   
{7} &  Regularity for $u_{a}= {\textstyle \frac {1}{2}}(u_{1}+u_{2})$ and %
                                                                 $v=\pm \ha (u_{1}-u_{2})$, Part II &{33} \cr   
{8} &  A frequency function for $v$ and the dimension  %
                                                                 of ${\mathcal K}_{u}$  &{36} \cr  
{9} &  Appendix: A simple connectivity lemma    &{42} \cr } }} 

\section*{Introduction}

In \cite{Wic08}, the second author established a local $C^{1, \alpha}$ partial regularity theory for
stable branched minimal hypersurfaces of multiplicity $< 3$. The main regularity theorem of
\cite{Wic08} in particular implies that locally near any point where it has one tangent cone equal to
a multiplicity 2 hyperplane, a stationary integral varifold arising as the weak limit of a sequence of
stable minimal hypersurfaces, each of which is immersed away from a closed set of singularities
(including branch points) of locally finite codimension 2 Hausdorff measure, must be a
$C^{1,\alpha}$ 2-valued graph over a domain in a suitable hyperplane for some $\alpha\in (0, 1)$. It
remained open how large the singular set of the varifold could be near such a point. Also left open
was the question of the optimal value of $\alpha$.

Here we give answers to these questions by proving that, near such a point, the varifold is always a
$C^{1,1/2}$ 2-valued graph, and that either it is regular (which means that in a neighborhood of the
point, either the support of the varifold decomposes as the union of two smooth
embedded (intersecting) graphs or the varifold is equal to a multiplicity 2 copy of
a single regular embedded minimal graph), or the set of its singularities (i.e.\ branch points) has
Hausdorff dimension precisely equal to $n-2$.

In fact we here establish, in Theorems~7.1,~7.4 and~8.10, that such results apply in arbitrary
codimension $k\ge 1$ without any {\it a priori\/} stability assumption. That is, we show that any
2-valued $C^{1, \alpha}$ ($\alpha \in (0, 1)$) function $u=\{u_{1},u_{2}\}$ (with values of
$u_{j}\in\R^{k}$) on an open ball $B$ in $\R^{n}$ whose graph $G$, viewed as a varifold with
multiplicity 2 at points where $u_{1}=u_{2}$ and with multiplicity 1 at points where $u_{1}\neq
u_{2}$, is stationary in the cylinder $B \times \R^{k}$, must be a $C^{1, 1/2}$ function, and that
the set of its singularities, if non-empty, must have Hausdorff dimension equal to $n-2$.  The
$C^{1, 1/2}$ regularity is of course optimal, as is shown by the case when $n=k=2$ and $u$ is the
2-valued function $u(x,y)=z^{3/2},\,z=x+iy$, with graph
$G=\{(w,z)\in\C\times\C\approx\R^{2}\times\R^{2}:w^{2}=z^{3}\}$ which is a complex algebraic
variety, hence minimizing in $\R^{4}$ as a $2$-dimensional multiplicity~1 current, and hence
stationary as a 2~dimensional multiplicity~1 varifold in~$\R^{4}$.

 The existence of large families of non-parametric branched $C^{1,\alpha}$ minimal hypersurfaces
(i.e.\ the case $n\ge 2$, $k=1$) has been established by the authors \cite{SimW05}. In the case
$n=2$, recently L.~Rosales \cite{Ros10} established the existence of further classes of such surfaces
without the symmetry assumptions needed in \cite{SimW05}.

The main tool used here to bound the size of the branch set is a monotone frequency function for
the 2-valued difference $v=\pm\ha (u_{1}-u_{2})$. The frequency function allows one to produce
non-trivial, homogeneous 2-valued stationary harmonic blow-ups at branch points. 
F.~J.~Almgren~Jr.\ first introduced the notion of frequency function in the 1970's and used it to
study energy minimizing multiple-valued harmonic functions and the singular set of area
minimizing currents.  Almgren's main work on these topics, available since the early 1980's in
preprint form, was published posthumously in book form in~\cite{Alm00}.

Establishing monotonicity properties of the frequency function in the present PDE setting depends
crucially on knowing the $C^{1, 1/2}$ regularity of the solution. In fact we need, and prove, more
than that. We show that the (single valued) average function $u_{a}=\fr 1/2 (u_{1}+u_{2})$ of the
2-valued solution is of class $C^{1, 1},$ and the 2-valued difference $v=\pm\ha (u_{1}-u_{2})$ is in
$C^{1, 1/2}$. Most of the present work goes into proving these regularity results, and involves in
particular establishing a $C^{1, \alpha}$ Schauder theory and $W^{2, 2}$ estimates for our 2-valued
functions, as well as a ``frequency gap'' result for 2-valued stationary harmonic functions and
growth estimates for 2-valued solutions to certain linear equations.  Once the required regularity is
established, it is straightforward to prove that the 2-valued difference function $v=\pm\ha
(u_{1}-u_{2})$ satisfies a weakly coupled divergence form elliptic system with Lipschitz coefficients,
and we can then apply appropriate modifications of the work of Garofolo and Lin~\cite{GarL86}
(which establishes monotonicity of a frequency function for single valued solutions of
divergence-form elliptic equations with Lipschitz coefficients).

The remainder of the proof depends on more or less standard application of ``dimension reducing''
arguments utilizing the monotonicity of the frequency function in a manner completely analogous
to the arguments in~\cite{Alm00}.

\section*{1 \quad Preliminaries}

We use the notation $B_{\rho}(x_{0})=\{x\in\R^{n}: |x-x_{0}|<\rho\}$,
$\ovl{B}_{\rho}(x_{0})=\{x\in\R^{n}: |x-x_{0}|\le \rho\}$, $k\ge 1$ ($k$ is the codimension),
$B_{\rho}=B_{\rho}(0)$, and $u$ denotes a $C^{1,\alpha}(B_{1},\R^{k})$ 2-valued function, so
$$%
u(x)=\{u_{1}(x),u_{2}(x)\} %
\text{ (an unordered pair of points in $\R^{k}$) for each $x\in B_{1}$. } %
\leqno{\bn{1.1}}
$$%
With such 2-valued functions $u=\{u_{1},u_{2}\},v=\{v_{1},v_{2}\}$ we adopt the convention that
\begin{align*}%
  |u(x)| %
  &= |u_{1}(x)| +|u_{2}(x)| \\ %
  |u(x)-v(x)| %
  &= \min\{|u_{1}(x)-v_{1}(x)|+|u_{2}(x)-v_{2}(x)|, %
  |u_{1}(x)-v_{2}(x)|+|u_{2}(x)-v_{1}(x)|\}%
\end{align*}%
and we write $u\in C^{1}(B_{1})$ if for each $x\in B_{1}$ there is a 2-valued affine function
$L_{x}$ on $\R^{n}$ of the form
$$%
L_{x}(h) =\{u_{1}(x)+A_{1}(x)h,u_{2}(x)+A_{2}(x)h\}
$$%
(assuming points in $\R^{n}$ are written as columns). Here $A_{1}(x),A_{2}(x)$ are $k\times n$
matrices and
$$%
\lim_{h\to 0}|h|^{-1}|u(x)- L_{x}(h)|=0,\,\, %
\lim_{y\to x}\sup_{|h|=1}|L_{x}(h)-L_{y}(h)|=0, %
$$%
and in this case $Du(x)$ denotes the (unique) 2-valued function $\{A_{1}(x),A_{2}(x)\}$, and we
sometimes also write $Du_{1}(x),Du_{2}(x)$ rather than $A_{1}(x),A_{2}(x)$. Also we say $u\in
C^{1,\alpha}(B_{1})$ if
$$%
|u|_{1,\alpha;B_{1}}<\infty, %
$$%
where $\alpha\in (0,1)$ is given and
$$%
|u|_{1,\alpha;B_{1}}=\sup_{B_{1}}|u|+\sup_{B_{1}}|Du| +[Du]_{\alpha,B_{1}},
$$%
where $|u|,|Du|$ and the H\"older coefficient $[Du]_{\alpha,B_{1}}$ are interpreted in the usual
way; thus
\al{%
  |u|&=|u_{1}|+|u_{2}|,\quad |Du|=|Du_{1}|+|Du_{2}| \\ %
  [Du]_{\alpha} &= \sup_{x_{1},x_{2}\in B_{1},\,x_{1}\neq x_{2}} %
  |x_{1}-x_{2}|^{-\alpha}|Du(x_{1})-Du(x_{2})|, }%
where
\begin{align*}%
  |Du(x_{1})-Du(x_{2})| %
  &=\min\{|A_{1}(x_{1})-A_{1}(x_{2})|+|A_{2}(x_{1})-A_{2}(x_{2})|,\\ %
  &\hskip2in |A_{1}(x_{1})-A_{2}(x_{2})|+|A_{2}(x_{1})-A_{1}(x_{2})|\}, %
\end{align*}%
with $A_{1}(x),A_{2}(x)$ as above.  The area functional is given by
$$%
{\cal{}A}(u) = \sint_{B_{1}}\Bl(\sqrt{g(u_{1})} + \sqrt{g(u_{2})}\Br)
$$%
where $g(u_{\ell})=\det(\delta_{ij}+D_{i}u_{\ell}\cdot D_{j}u_{\ell})$ (notice this makes sense
because $\sqrt{g(u_{1})} + \sqrt{g(u_{2})}$ is a well defined single-valued function on $B_{1}$), and
we assume that $u$ is a stationary point for this functional in the sense that $G=\graph u=\{(x,y)\in
B_{1}\times\R^{k}:y=u_{1}(x)\text{ or }y=u_{2}(x)\}$ is a stationary varifold. Thus we assume that
$$%
\sint_{G} \dvg_{G\!} X \,\theta d{\cal{}H}^{n}= 0,\quad j=1,\ldots,n+k, %
X\in C^{1}_{\text{c}}(B_{1}\times \R^{k},\R^{n+k}), %
\leqno{\bn{1.2}}
$$%
where $\theta$ is the multiplicity function ($=2$ and points where $u_{1}=u_{2}$ and $=1$ at
points where $u_{1}\neq u_{2}$) and where $\dvg_{G\!}X$ denotes the tangential divergence of $X$
on $G$. Thus $\dvg_{G}X=\sum_{j=1}^{n+k}e_{j}\cdot \nabla_{G}X_{j}$, with $\nabla_{G}X_{j}$
denoting the gradient of $X_{j}$ on $G$ (i.e.\ $P_{x}(D X_{j})$, where $P_{x}$ is the orthogonal
projection of $\R^{n+k}$ onto the tangent space of $G$ at any point $x\in G$). In particular if we
let
$$%
{\cal{}K}_{u}=\{x\in B_{1}: u_{1}(x)=u_{2}(x)\text{ and }Du_{1}(x)=Du_{2}(x)\} %
$$%
then in each ball $B_{\sigma}(y)\subset B_{1}\setminus {\cal{}K}_{u}$ we can label the values
$u_{1},u_{2}$ of $u$ such that $u_{1}|B_{\sigma}(y),u_{2}|B_{\sigma}(y)$ are
$C^{\infty}(B_{\sigma}(y))$ solutions of the minimal surface system, so that
$$%
{\cal{}M}(u_{j}) =0 \text{ in } B_{\sigma}(y),\quad j=1,2, %
\leqno{\bn{1.3}}
$$%
with ${\cal{}M}(u_{j})=({\cal{}M}_{1}(u_{j}),\ldots,{\cal{}M}_{k}(u_{j}))$,
$$%
{\cal{}M}_{\kappa}(w)= %
\ssum_{j=1}^{n}D_{j}\bigl(\sqrt{g(w)}g^{ij}(w)D_{j}w_{\kappa}\,\bigr),%
\quad \kappa=1,\ldots,k, %
\leqno{\bn{1.4}}
$$%
where $(g^{ij}(w))=(g_{ij}(w))^{-1}$, $g_{ij}(w)=\delta_{ij}+D_{i}w\cdot D_{j}w$, $g(w)=\det
(g_{ij}(w))$. Our aim is to show that the closed set ${\cal{}K}_{u}$ in fact has Hausdorff dimension
$\le n-2$. Observe that, since $u$ is $C^{1}$ and 2-valued one can check the inclusion
$$%
{\cal{}B}_{u}\subset {\cal{}K}_{u}, %
\leqno{\bn{1.5}}
$$%
where ${\cal{}B}_{u}$ is the ``branch set'' ${\cal{}B}_{u}$ of $u$, which is defined to be the set of
points $y\in B_{1}$ such that there is no neighborhood $U_{y}$ of $y$ such that the values
$u_{1},u_{2}$ can be ordered in $U_{y}$ in such a way that each of $u_{1},u_{2}$ is a single-valued
$C^{1}$ function in $U_{y}$. Observe that the inclusion~1.5 fails in general for $Q$-valued
$C^{1}$ functions with $Q\ge 3$.

The function
$$%
v(x)=\{\pm \ha (u_{2}(x)-u_{1}(x)):x\in B_{1}\}
$$%
defines a 2-valued $C^{1,\alpha}(B_{1},\R^{k})$ function which is ``symmetric,'' in the sense that at
each point $x$ the two values of $v(x)$ are negatives of each other.  ${\cal{}K}_{u}$ is then the same
as
$$%
{\cal{}K}_{v}=\{x\in B_{1}: |v(x)|=0,\,|Dv(x)|=0\}.
$$%
The main results proved here are local results valid in a neighborhood of a point ($0$ say) in
${\cal{}K}_{u}$ and since we can always (rotating the graph if necessary) assume that
$u(0)=\{0,0\}, Du(0)=\{0,0\}$, we can, and we shall, assume (after rescaling) that in fact
$$%
\sup_{B_{1}}|u| + \sup_{B_{1}}|Du| + [Du]_{\alpha,B_{1}} \le \epsilon_{0}, \leqno{\bn{1.6}}
$$%
where $\epsilon_{0}$ is to be specified (depending only on $n,k$) later.  Notice that by using~1.6
together with (single-valued) quasilinear elliptic estimates in balls contained $B_{1}\setminus
{\cal{}K}_{v}$, we have for each $\theta\in (0,1)$
$$%
|v(x)|+d(x)|Dv(x)|+d(x)^{2}|D^{2}u(x)| \le C \epsilon_{0}d(x)^{1+\alpha} %
\leqno{\bn{1.7}}
$$%
on $B_{\theta}$, where, here and subsequently,
$d(x)=\dist(x,{\cal{}K}_{v})\,(=\dist(x,{\cal{}K}_{u}))$, and where $C=C(n,\theta)$.  There is a
well-defined single-valued $C^{1,\alpha}$ ``average'' $u_{a}$ given by
$$%
u_{a}=\ha(u_{1}+u_{2}) \text{ in } B_{1}, %
\leqno{\bn{1.8}}
$$%
where $u_{1}, u_{2}$ are as in~1.1.  A principal ingredient in the proof that the frequency function
for $v$ has the appropriate monotonicity properties will involve showing that $u_{a}$ is of class
$C^{1,1}$ (which we do in \S7).

As mentioned above, the function $v=\{\pm\ha (u_{1}-u_{2})\}$ is a 2-valued $C^{1,\alpha}$
symmetric function, and for this reason much of the analysis that follows will relate to 2-valued
symmetric functions. Since we use integral estimates it is necessary to discuss Sobolev spaces of such
functions. So assume that $w$ is a 2-valued symmetric (i.e.\ at each point the two values of $w$ are
negatives of each other), ${\cal{}Z}_{w}=\{x\in B_{1}:w(x)=\{0,0\}\}$, $w\in C^{0}(B_{1})\cap
C^{1}(B_{1}\setminus {\cal{}Z}_{w})$, and observe that in any ball contained in $B_{1}\setminus
{\cal{}Z}_{w}$, we can represent $w$ uniquely as $\pm w_{1}$ for some unique positive $C^{0}$
function $w_{1}$. We say that $w\in W^{1,p}(B_{1})$ if $D_{j}w\in L^{p}(B_{1})$, where
$D_{j}w$ is the symmetric 2-valued function defined locally near a point $\xi\in B_{1}\setminus
{\cal{}Z}_{w}$ as $\pm D_{j}w_{1}$ and $D_{j}w$ is defined to be $\{0,0\}$ on ${\cal{}Z}_{w}$.

In practice it is usually more convenient to use the equivalent definition
$$%
D_{j}w=\lim_{\delta\downarrow 0} D_{j}\gamma_{\delta}(w) %
\leqno{\bn{1.9}}
$$%
(limit taken in $L^{p}$), where $\delta>0$ and $\gamma_{\delta}$ denotes a smooth odd
($\gamma_{\delta}(-t)=-\gamma_{\delta}(t)$) increasing function on $\R$ with the properties that
$\gamma_{\delta}$ vanishes identically in some neighborhood of $0$,
$\gamma_{\delta}^{\prime}(t)\le 1$ for all $t$, $\gamma_{\delta}(t)\equiv t-\delta$ for $t>\delta$,
$\gamma_{\delta}(t)=t+\delta$ for $t<-\delta$. Using this characterization one can easily check for
example that then $|w|^{2}\in W^{1,p}(B_{1})$ (as a single-valued function) with weak derivatives
$2w\cdot D_{j}w$, assuming that we adopt the natural convention that, near points $\xi\in
B_{1}\setminus {\cal{}Z}_{w}$ where we can in a unique way write $w=\pm w_{1}$ with $w_{1}$
continuous,
$$%
\text{$|w|^{2}$ and $w\cdot D_{j}w$ are taken to mean $|w_{1}|^{2}$ and %
  $w_{1}\cdot D_{j}w_{1}$ respectively} %
\leqno{\bn{1.10}}
$$%
on $B_{1}\setminus {\cal{}Z}_{w}$ (and $=0$ on $B_{1}\cap {\cal{}Z}_{w}$).  Note that if $w\in
C^{1}(B_{1})\cap C^{2}(B_{1}\setminus {\cal{}Z}_{w})$ is symmetric then the (classical or weak)
derivatives $D_{j}w$ are in $C^{0}(B_{1})\cap C^{1}(B_{1}\setminus {\cal{}Z}_{w})$ and it therefore
makes sense to define the second order weak derivative $D_{i}D_{j}w$ by $D_{i}(D_{j}w)$ in
accordance with the above discussion with $D_{j}w$ in place of $w$.  One then easily checks (using
approximation involving $\gamma_{\delta}$ as above) that for example if $D_{i}D_{j}w\in
C^{1}(B_{1}\setminus {\cal{}Z}_{w})$ then $D_{\ell}(D_{i}wD_{j}w)=D_{\ell}D_{i}w \, D_{j}w+
D_{i}w\,D_{\ell}D_{j}w$ on $B_{1}\setminus{\cal{}K}_{w}$, assuming that we define the products
naturally as in~1.10 above.

\section*{2\quad 2-valued {\boldmath $C^{1}$} harmonic %
  functions---Part I}

Given 2-valued symmetric $C^{1}(B_{1},\R^{k})$ function $\varphi$, we say that $\varphi$ is
harmonic if for each ball $B_{\sigma}(y)\subset B_{1}\setminus{\cal{}K}_{\varphi}$ there is $C^{1}$
harmonic function $\varphi_{1}$ on $B_{\sigma}(y)$ such that $\varphi|B_{\sigma}(y)=\{\pm
\varphi_{1}\}$; of course if it exists such a $\varphi_{1}$ is unique.

Our first aim is to show that such symmetric $C^{1}$ functions are automatically locally $W^{2,2}$
in $B_{1}$, with an estimate on the $W^{2,2}$ norm and the Lipschitz constant.

\begin{state}{\bf{}2.1 Lemma.} %
  Suppose that $\varphi$ is a $C^{1}$ \emph{2}-valued symmetric harmonic function on $B_{1}$ (in
the above sense).  Then $D^{2}\varphi\in L^{2}(B_{\rho}(y))$ (i.e.\ $D\varphi\in
W^{1,2}(B_{\rho}(y))$) for each ball $B_{\rho}(y)$ with $\ovl{B}_{\rho}(y)\subset B_{1}$, and we
have the estimates
$$%
\rho^{2-n}\sint_{B_{\rho/2}(y)}|D^{2}\varphi|^{2} + %
\sup_{B_{\rho/2}(y)}|D\varphi|^{2} \le C\rho^{-n}\sint_{B_{\rho}(y)}|D\varphi|^{2}
$$%
for all such balls $B_{\rho}(y)$.
\end{state}%

{\bf{}Proof:} Since $\varphi=\{\pm \varphi_{1}\} $ in any ball away from ${\cal{}K}_{\varphi}$, we
have $\Delta |D\varphi|^{2}=2|D^{2}\varphi|^{2}$ on $B_{1}\setminus{\cal{}K}_{\varphi}$, so if
$\gamma_{\delta}$ is a smooth non-negative convex function on $\R$ with $\gamma_{\delta}\equiv
0$ in some neighborhood of zero and $\gamma_{\delta}^{\prime}(t)\equiv 1$ on $[\delta,\infty)$ for
some $\delta>0$, then we have that
$$%
\Delta \gamma_{\delta}(|D\varphi|^{2})\ge
2\gamma_{\delta}^{\prime}(|D\varphi|^{2})|D^{2}\varphi|^{2} %
$$%
on any closed ball $B\subset B_{1}$ (because $\gamma_{\delta}(|D\varphi|^{2})|B$ has compact
support in $B\setminus{\cal{}K}_{\varphi}$), and so multiplying by a cut-off function which is
identically 1 in $B_{\rho/2}(y)$ and zero outside $B_{\rho}(y)$, and integrating over $B_{\rho}(y)$,
we obtain the required $W^{2,2}$ estimate by letting $\delta\downarrow 0$.  Also the above
inequality shows that $\gamma_{\delta}(|D\varphi|^{2})$ is a subharmonic (single-valued) function
in $B_{1}$, and so for each ball $B_{\rho}(y)\subset B_{1}$ we have the estimate
$\sup_{B_{\rho/2}(y)}\gamma_{\delta}(|D\varphi|^{2})\le
C\rho^{-n}\int_{B_{\rho}(y)}\gamma_{\delta}(|D\varphi|^{2})$ and, again we can let
$\delta\downarrow 0$ to get the required estimate for $\sup_{B_{\rho/2}(y)}|D\varphi|^{2}$.

\medskip

\begin{state}{\bf{}2.2 Lemma.} %
  Suppose that $\varphi$ is a $C^{1}$ \emph{2}-valued symmetric function on the ball
$B_{R}(y)\subset\R^{n}$, that $\varphi$ is harmonic on $B_{R}(y)\setminus {\cal{}K}_{\varphi}$ and
not identically zero on $B_{R}(y)$, and that $y\in{\cal{}Z}_{\varphi}=\{x:\varphi(x)=\{0,0\} \}$.  Then
$$%
N_{\varphi}(y,\rho)=\Fr\rho^{2-n}\mint_{B_{\rho}(y)}|D\varphi|^{2}/ %
{\rho^{1-n}\mint_{\partial B_{\rho}(y)}|\varphi|^{2}}
$$%
is an increasing $C^{1}$ function of $\rho\in (0,R)$, ${\cal{}N}_{\! 
\varphi}(y)=\lim_{\rho\downarrow 0}N_{\varphi}(y,\rho)\ge 1$, and ${\cal{}Z}_{\varphi}$ (and hence
${\cal{}K}_{\varphi}$) has empty interior.
\end{state}%

{\bf{} 2.3 Remarks: (1)} $N_{\varphi}(y,\rho)$ is called the frequency function of $\varphi$,
terminology introduced by Almgren~\cite{Alm00}. The frequency function was a key tool in
Almgren's study of energy minimizing multi-valued functions and area minimizing currents. As
observed by Almgren, the fact that $N_{\varphi}(y,\rho)$ is increasing is equivalent to the fact that
$$%
\log\Bigl(\rho^{1-n}\sint_{\partial B_{\rho}(y)}|\varphi|^{2}\Bigr) \text{ is a convex function of
}\,t=\log \rho,
$$%
because (as discussed in~(2) below)
$$%
\rho^{2-n}\sint_{B_{\rho}(y)}|D\varphi|^{2}= \fr 1/2 %
\Fr d/{d\rho}\Bigl(\rho^{1-n}\sint_{\partial B_{\rho}(y)}|\varphi|^{2}\Bigr), %
$$%
and hence $N_{\varphi}(y,\rho)$ can be written alternatively
$$%
N_{\varphi}(y,\rho) = \frac{\fr 1/2 \rho H^{\prime}(y,\rho)}{H(y,\rho)}, %
\quad H(y,\rho)=\rho^{1-n}\int_{\partial B_{\rho}(y)}|\varphi|^{2}.
$$%

{\bf{}(2)} The proof of the monotonicity of $N_{\varphi}$ will be based on the key identities
\begin{align*}%
  \sint_{B_{\rho}(y)}|D\varphi|^{2} %
  &= \sint_{\partial B_{\rho}(y)}\varphi\cdot D_{r}\varphi \,\,\,\, %
  \Bigl(= \fr 1/2\sint_{\partial B_{\rho}(y)} D_{r}|\varphi|^{2} \Bigr) \\ %
  \sint_{B_{\rho}(y)}\Bigl(|D\varphi|^{2}\delta_{ij} - %
  2D_{i}\varphi\cdot D_{j}\varphi\Bigr) D_{i}\zeta_{j} %
  &=\sint_{\partial B_{\rho}(y)}\Bigl(|D\varphi|^{2}(\rho^{-1}(x-y)\cdot \zeta) - %
  2D_{r}\varphi\cdot D_{j}\varphi\zeta_{j}\Bigr) %
\end{align*}%
where $\zeta_{j}$ are Lipschitz in $\ovl{B}_{\rho}(y)$, where $\ovl B_{\rho}(y)\subset B_{1}$ and
where the radial derivative $D_{r}\varphi_{\ell}$ is defined by $\rho^{-1}(x-y)\cdot D\varphi_{\ell}$
for $\ell=1,\ldots,k$.  These identities are readily checked using integration by parts, which is
justified by virtue of Lemma~2.1.

{\bf{}(3)} In view of the monotonicity in Lemma~2.2, it follows directly that ${\cal{}N}_{\! 
\varphi}(y)$ is an upper semicontinous function of $\varphi$ with respect to the $W^{1,2}$-norm;
thus if $\varphi_{j}\to \varphi$ in $W^{1,2}(B_{\rho}(y))$ then ${\cal{}N}_{\! \varphi}(y)\ge
\limsup_{j\to \infty} {\cal{}N}_{\! \varphi_{j}}(y)$. By applying this to the functions
$\varphi_{j}(x)=\varphi(x+y_{j}-y)$ we also have ${\cal{}N}_{\! \varphi}(y)\ge \limsup_{y_{j}\to
y}{\cal{}N}_{\! \varphi}(y_{j})$.

{\bf{}(4)} Note also that, with $H(y,\rho)$ as in~(1) above, we have the following general growth
facts related to $N$ for $\varphi$ as in~2.2: {\abovedisplayskip3pt\belowdisplayskip1pt
$$%
{\cal{}N}_{\varphi}(y)\le \tfrac{\rho}{2} H^{\prime}(y,\rho)/{H(y,\rho)}=N_{\varphi}(y,\rho)\le
N_{\varphi}(y,R), \quad \rho\in (0,R],
$$}%
and hence by integration we have the bounds
$$%
\Bigl(\Fr \rho/R\Bigr)^{C} \le \sqrt{\Fr H(y,\rho)/{H(y,R)}} \le %
\Bigl(\Fr \rho/R\Bigr)^{{\cal{}N}_{\varphi}(y)}, \quad \rho\in (0,R], \, C=N_{\varphi}(y,R).
$$%
{\bf{}(5)} Notice that~(4) above (with $R=2\rho$) and the monotonicity of $N_{\varphi}(y,\sigma)\le
N_{\varphi}(y,R)$ for $\sigma\in (0,R]$ imply that $C^{-1}H(y,2\rho) \le H(y,\rho) \le C
H(y,2\rho)$ for $\rho\in (0,R/2]$ with fixed $C=2^{N_{\varphi}(y,R))}$ and so by integrating with
respect to $\rho$ we have
$$%
\|\varphi \|_{L^{2}(B_{2\rho}(y))}\le C\|\varphi\|_{L^{2}(B_{\rho}(y))},\quad \rho\in (0,R/2],
$$%
again with fixed $C=C(N_{\varphi}(y,R))$.

\medskip

{\bf{}Proof of 2.2:} Observe first that by the first identity in~2.3(2) we see that $|\varphi|\equiv 0$ on
$\partial B_{\rho}(y)\Rightarrow \varphi|B_{\rho}(y)\equiv 0$, so the frequency
$N_{\varphi}(z,\rho)$ is well-defined unless $\varphi\equiv 0$ in $B_{\rho}(y)$. Also, taking
$\zeta_{j}\equiv x_{j}-y_{j}$ in the second identity of~2.3(2), we see that
$$%
{d\over{}d\rho}\Bigl(\rho^{2-n}\int_{B_{\rho}(y)}|D\varphi|^{2}\Bigr) = %
2\rho^{2-n}\int_{\partial B_{\rho}(y)}|D_{r}\varphi|^{2} %
\Bigl(= 2\rho^{}\int_{S^{n-1}}|D_{r}\varphi(y+r\omega)|^{2}\,d\omega\Bigr); %
\leqno{(1)}
$$%
Notice that this can be alternatively written
$$%
{dD(y,\rho)\over{}d\rho}=2\rho\int _{S^{n-1}} |D_{r}\varphi(y+r\omega)|^{2}\,d\omega,
$$%
with $D(y,\rho)=\rho^{2-n}\int_{B_{\rho}(y)}|D\varphi|^{2}$, and the first identity in~2.3(2) can be
written
$$%
D(y,\rho)=\ha{dH(y,\rho)\over{}d\rho}
$$%
Assuming $\rho_{0}\in (0,R)$ with $\varphi|B_{\rho_{0}}(y)$ not identically zero, one can now
directly check from these two identities that
$$%
{\fr d/{d\rho}}N(y,\rho)= %
2H(y,\rho)^{-2}\rho\Bigl(\int_{S^{n-1}}\varphi_{r}^{2}\int_{S^{n-1}}|\varphi|^{2}- %
\Bigl(\int_{S^{n-1}}\varphi\varphi_{r}\Bigr)^{2}\Bigr) ,\quad \rho\in (\rho_{0},R), %
$$%
and, since the right side here is non-negative by the Cauchy-Schwarz inequality, we thus have the
monotonicity
$$%
{\fr d/{d\rho}}N(y,\rho) \ge 0 ,\quad \rho\in (\rho_{0},R).  %
\leqno{(2)}
$$%
Furthermore it then follows that
$$%
\text{${\cal{}Z}_{\varphi}$ has empty interior}, %
\leqno{(3)}
$$%
because otherwise, since $\varphi$ is not identically zero by hypothesis, we could find
$\sigma,\delta>0$ and balls $B_{\sigma}(z)\subset B_{\sigma+\delta}(z)\subset B_{R}(y)$ with
$|\varphi||B_{\sigma}(z)\equiv 0$ but $\sup_{B_{\rho}(z)}\varphi|B_{\rho}(z)>0$ for all
$\rho\in(\sigma,\sigma+\delta)$, and (see the discussion in Remark~2.4(1) above) since ${\fr
d/{d\rho}}N(y,\rho)\ge 0$ can be written ${\fr d^{2}/{dt^{2}}}\log (\rho^{1-n}\int_{\partial
B_{\rho}}\varphi^{2})\ge 0$, where $t=\log \rho$, we see that $\log (\rho^{1-n}\int_{\partial
B_{\rho}}\varphi^{2})$ is bounded below as $\rho\downarrow\sigma$ (because any convex function
on an open interval is bounded below), contradicting the fact that $\rho^{1-n}\int_{\partial
B_{\rho}}\varphi^{2}\to 0$ as $\rho\downarrow \sigma$.  Hence we actually have
$$%
\text{$N_{\varphi}(y,\rho)$ is a well-defined %
  $C^{1}$ increasing function of $\rho$ for $\rho\in (0,R)$} %
\leqno{(4)}
$$%
as claimed.

Finally, if ${\cal{}N}_{\varphi}(y)<1$ then we could choose $\rho_{0}\in (0,R)$ such that
$N_{\varphi}(y,\rho_{0})=C_{0}<1$ and then Remark~2.3(4) would give $\sqrt{H(\rho)}\ge
C\rho^{C_{0}}$ as $\rho\downarrow 0$, whereas since $\varphi$ is $C^{1}$ with
$\varphi(0)=\{0,0\}$ we must have $\sqrt{H(\rho)}\le C\rho$ as $\rho\downarrow 0$. Thus we also
have ${\cal{}N}_{\!  \varphi}(y)\ge 1$ for each $y\in {\cal{}Z}_{\varphi}$ as claimed.

\medskip

{\textbf{2.4 Remarks. (1)\,}} Note that by examining the proof of the monotonicity of
$N_{\varphi}(y,\rho)$ in the first part of the above proof, we see that $N_{\varphi}(y,\rho)$ can be
constant in some interval $(\rho_{0},\rho_{0}+\epsilon)$ ($\epsilon>0$) if and only if we have
equality in the Cauchy-Schwarz inequality $\int_{S^{n-1}}\varphi_{r}^{2}\int_{S^{n-1}}|\varphi|^{2}-
\Bigl(\int_{S^{n-1}}\varphi\varphi_{r}\Bigr)^{2}\ge 0$ which in turn is true if and only if
$\varphi_{r}$ is a constant multiple of $\varphi$ for $|x-y|\in (\rho_{0},\rho_{0}+\epsilon)$, which is
in turn true if and only if $\varphi$ extends to be a homogeneous function with respect to the
variable $r=|x-y|$, with degree $\beta=$ the constant value of $N_{\varphi}(y,\rho)$.

{\bf{}(2)\,} In view of Remark~(1) above we see that if $\varphi$ is a 2-valued $C^{1}$ symmetric
function which is harmonic on $\R^{n}\setminus {\cal{}K}_{\varphi}$ and homogeneous of degree
$\beta\ge 1$ (i.e.\ $\varphi(\lambda x)=\lambda^{\beta}\varphi(x)$), then ${\cal{}N}_{\varphi}(z)\le
{\cal{}N}_{\varphi}(0)$ for each $z\in\R^{n}$ and $S=\{x\in
{\cal{}Z}_{\varphi}:{\cal{}N}_{\varphi}(z)={\cal{}N}_{\varphi}(0)\}$ is a linear subspace with
$\varphi\circ \tau_{z}=\varphi$ for each $z\in S$, where $\tau_{z}$ is the translation $x\mapsto x+z$.
 This does follow directly from~Remark~(1) if one keeps in mind that if we have $z\neq 0$ and the
two homogeneity conditions $\varphi(\lambda x)=\lambda^{\beta}\varphi(x)$ and $\varphi(\lambda
x + z)=\lambda^{\beta}(x+z)$ for each $\lambda>0$ and each $x\in\R^{n}$, then we have, with $t\in
\R$ arbitrary and $\lambda>0$ chosen so that $\lambda^{-1}-\lambda=t$,
$\varphi(x+tz)=\varphi(x-\lambda z+ \lambda^{-1}z)= \lambda^{-\beta}\varphi(\lambda
x-\lambda^{2}z+z)=\lambda^{\beta}\varphi(\lambda^{-1}x-z+z)=
\lambda^{\beta}\varphi(\lambda^{-1}x)=\varphi(x)$.

\medskip

The following gap lemma will be used to establish a Liouville-type theorem for symmetric 2-valued
harmonic functions (in 2.6 below), which in turn will be the main ingredient in the $C^{1,\alpha}$
Schauder theory of~\S3.

\begin{state}{\bf{}2.5 Lemma.} %
  There is $\delta=\delta(n)\in (0,1)$ such that if $\varphi(x)=|x|^{\sigma}\varphi(|x|^{-1}x)$ is
2-valued symmetric $C^{1}$ homogeneous degree $\sigma$ function with $\sigma\in [1,1 +\delta)$
and with $\varphi$ harmonic on $\R^{n}\setminus {\cal{}K}_{\varphi}$ and not identically zero,
then $\sigma=1$ and $\varphi$ is linear (i.e.\ $\varphi(x)\equiv\{\pm \ell(x)\}$, where
$\ell_{\kappa}(x)=\sum_{j=1}^{n}c_{\kappa}^{j}x_{j}$ for some constants $c_{\kappa}^{j}$,
$\kappa=1,\ldots,k$).
\end{state}%

{\bf{}Proof:} We can assume that $\varphi$ is real-valued (i.e.\ $k=1$) because each component
$\varphi_{\kappa}$ is either identically zero or satisfies the stated hypotheses with $k=1$ and
$\varphi_{\kappa}$ in place of $\varphi$.

The theorem is trivially true in case $n=1$, so assume $n\ge 2$.  We first dispense with the case
$\sigma=1$. For $\delta>0$ let $\gamma_{\delta}(t)$ be an odd function of $t$ which is convex for
$t\ge 0$, $\equiv 0$ in some neighborhood of $0$, and which has
$\gamma_{\delta}^{\prime}(t)\equiv 1$ for $t\ge\delta$, and observe that
$\gamma_{\delta}(D_{j}\varphi)$ has compact support in $S^{n-1}\setminus{\cal{}K}_{\varphi}$,
and hence (since $D_{j}\varphi$ is harmonic and homogeneous degree $0$ on $\R^{n}\setminus
{\cal{}K}_{\varphi}$) we have, interpreting products in the natural way on $\R^{n}\setminus
{\cal{}K}_{\varphi}$ (Cf.~1.10),
$$%
\sint_{S^{n-1}}|\nabla_{S^{n-1}}\gamma_{\delta}(D_{j}\varphi)|^{2} \le %
\sint_{S^{n-1}}\nabla_{S^{n-1}}D_{j}\varphi\cdot %
\nabla_{S^{n-1}}\gamma_{\delta}(D_{j}\varphi) = %
-\sint_{S^{n-1}}(\Delta_{S^{n-1}}D_{j}\varphi)\, \gamma_{\delta}(D_{j}\varphi)=0,
$$%
and hence $\gamma_{\delta}(D_{j}\varphi)$ is $\{\pm c_{j}(\delta)\}$ for some constant
$c_{j}(\delta)$, and so $D_{j}\varphi$ is $\{\pm c_{j}\}$ for some constant $c_{j}$ and the lemma is
proved in case $\sigma=1$.

To prove the case $\sigma\in (1,1+\delta)$, observe if there is no $\delta$ as claimed, then we would
have a sequence $\varphi^{(j)}$ of 2-valued symmetric $C^{1}$ functions, harmonic in
$\R^{n}\setminus{\cal{}K}_{\varphi^{(j)}}$, not identically zero, and homogeneous of degree
$\sigma_{j}$ with $\sigma_{j}>1$ and $\sigma_{j}\downarrow 1$. Let
${\cal{}K}_{j}={\cal{}K}_{\varphi^{(j)}}$. Assume without loss of generality that we have normalized
$\varphi^{(j)}$ so that $\|\varphi^{(j)}\|_{L^{2}(B_{1})}=1$ for each $j$. Then by~2.1 we have local
$W^{2,2}$ and Lipschitz estimates for $\varphi^{(j)}$ and a subsequence (still denoted
$\varphi^{(j)}$) converges locally uniformly and locally weakly in $W^{2,2}(\R^{n})$ to 2-valued
symmetric Lipschitz $W^{2,2}$ homogeneous degree 1 function $\varphi$.

With $\gamma_{\delta}$ as in the first part above, using the fact that $D_{i}\varphi^{(j)}$ is
harmonic and homogeneous degree $\sigma_{j}-1$ on $\R^{n}\setminus {\cal{}K}_{j}$ we have
\begin{align*}%
  \tag*{(1)} %
  \sint_{S^{n-1}}|\nabla_{S^{n-1}}\gamma_{\delta}(D_{i}\varphi^{(j)})|^{2} &\le %
  \sint_{S^{n-1}}\nabla_{S^{n-1}}D_{i}\varphi^{(j)}\cdot %
  \nabla_{S^{n-1}}\gamma_{\delta}(D_{i}\varphi^{(j)})\\ %
  &= (\sigma_{j}-1)(\sigma_{j}+n-3)\sint_{S^{n-1}} %
  D_{i}\varphi^{(j)}\gamma_{\delta}(D_{i}\varphi^{(j)}) \\ %
  &\le(\sigma_{j}-1)(\sigma_{j}+n-3)\sint_{S^{n-1}} |D_{i}\varphi^{(j)}|^{2}\to 0 %
  \text{ as }k\to \infty, %
\end{align*}%
and so, by Rellich's theorem, in the limit as $j\to \infty$ we conclude
$\nabla_{S^{n-1}}\gamma_{\delta}(D_{i}\varphi)=0$ a.e.\ for each $\delta>0$, so in fact
$D_{i}\varphi$ is given by $\pm c_{i}$ for some constant $c_{1},\ldots,c_{n}$. Modulo composition
with an orthogonal transformation of $\R^{n}$ we thus have
$$%
\text{$\varphi(x)\equiv\{\pm c_{0}x_{1}\}$ for some non-zero constant $c_{0}$.}
$$%
Observe that for each $\sigma, R>0$ we have $B_{R}\cap {\cal{}K}_{j}\cap
\R^{n}_{\sigma}=\emptyset$, and in fact $B_{R}\cap{\cal{}Z}_{j}\cap \R^{n}_{\sigma}=\emptyset$,
for all sufficiently large $j$, where $\R^{n}_{\sigma}=\{x\in\R^{n}:|x_{1}|> \sigma\}$ and
${\cal{}Z}_{j}=\{x:\varphi^{(j)}(x)=\{0,0\}\}$, and
$$%
D_{1}\varphi^{(j)}\to \{\pm c_{0}\} \text{ in }\R^{n}_{\sigma}, \quad %
D_{j}\varphi^{(j)}\to \{0,0\} \text{ in }\R^{n}_{\sigma} \text{ for $j=2,\ldots,n$}.  %
\leqno{(2)}
$$%
Observe also that if there is a constant $c_{j}$ such that
$D_{2}\varphi^{(j)}-c_{j}D_{1}\varphi^{(j)}=\{0,0\}$\hskip0.5pt\footnote[3]{Note that
$D_{2}\varphi^{(j)}-c D_{1}\varphi^{(j)}$ is to be interpreted in the natural way on
$\R^{n}\setminus {\cal{}K}_{j}$ as $\pm(D_{2}\varphi_{1}-c D_{1}\varphi_{1})$, where as usual we
locally write, on $\R^{n}\setminus {\cal{}K}_{j}$, $\varphi^{(j)}=\pm\varphi_{1}$ with $\varphi_{1}$
in $C^{1}$} on $\R^{n}\setminus{\cal{}K}_{j}$ (hence on all of $\R^{n}$), then $\varphi^{(j)}$ is
invariant under composition with translations in the direction $(-c_{j},1,0)$ and so for $n\ge 2$ we
could reduce the proof of the theorem from dimension $n$ to $n-1$. Therefore, since the theorem
is trivially true in case $n=1$, we can henceforth assume without loss of generality that for all $j$,
$$%
\text{ $\nexists\ c_{j}\in\R$ such that %
  $D_{2}\varphi^{(j)}-c_{j}D_{1}\varphi^{(j)}=\{0,0\}$ on %
  $\R^{n}\setminus{\cal{}K}_{j}$.} %
\leqno{(3)}
$$%
Now, with $\sigma_{0}\in (0,\fr 1/2)$ fixed, select a constant $c_{j}(\to 0)$ such that
$$%
\sint_{S^{n-1}\cap R^{n}_{\sigma_{0}}}(D_{2}\varphi^{(j)}- %
c_{j}D_{1}\varphi^{(j)}) D_{1}\varphi^{(j)}=0, %
\leqno{(4)}
$$%
and define
$$%
\psi^{(j)}=\|D_{2}\varphi^{(j)}- %
c_{j}D_{1}\varphi^{(j)}\|_{L^{2}(S^{n-1})}^{-1} (D_{2}\varphi^{(j)}-c_{j}D_{1}\varphi^{(j)}).
$$%
This is well-defined for all sufficiently large $j$ by~(3).  Since $\psi^{(j)}$ is harmonic and
homogeneous of degree $\sigma_{j}-1$ on $\R^{n}\setminus{\cal{}K}_{j}$, we can use the argument
of~(1), with $\psi^{(j)}$ in place of $D_{j}\varphi^{(j)}$, in order to conclude
\begin{align*}%
  \tag*{(5)} \sint_{S^{n-1}}|\nabla_{S^{n-1}}\gamma_{\delta}(\psi^{(j)})|^{2} \le %
  &\sint_{S^{n-1}}\nabla_{S^{n-1}}\psi^{(j)}\cdot %
  \nabla_{S^{n-1}}\gamma_{\delta}(\psi^{(j)})= \\ %
  &(\sigma_{j}-1)(\sigma_{j}+n-3)\sint_{S^{n-1}} \psi^{(j)} %
  \gamma_{\delta}(\psi^{(j)})\to 0\text{ as }j\to \infty
\end{align*}%
for each $\delta>0$ and hence $\psi^{(j)}$ converges locally weakly in $W^{1,2}$ to the symmetric
$\psi$ with $\psi(x)\equiv\{\pm c\}$) where $c$ (constant) is not zero because
$\int_{S^{n-1}}|\psi|^{2}=1$. However by multiplying by $\|D_{2}\varphi^{(j)}-
c_{j}D_{1}\varphi^{(j)}\|_{L^{2}(S^{n-1})}^{-1}$ in~(4), taking the limit in $j$ (keeping in mind~(2))
we then have $cc_{0}{\cal{}H}^{n-1}(S^{n-1}\cap \R^{n}_{\sigma_{0}})=0$, contradicting the fact
that $c$ and $c_{0}$ are non-zero constants.

\medskip

As a consequence of the Lemma~2.5 and Remark~2.3(1) we have the following Liouville-type result
for 2-valued symmetric $C^{1,\alpha}$ harmonic functions:

\begin{state}{\bf{}2.6 Corollary.} %
  If $\delta=\delta(n)$ is as in~2.5 and $\varphi$ is a 2-valued symmetric $C^{1,\alpha}$ function
on\/ $\R^{n}$ with $\alpha\in (0,\delta)$, $\varphi$ harmonic on\/
$\R^{n}\setminus{\cal{}K}_{\varphi}$ and $[D\varphi]_{\alpha,\R^{n}}<\infty$, then $\varphi$ is
affine (i.e.\ $\pm \ell$, where $\ell(x)$ is an affine function $c_{0}+\sum_{j=1}^{n}c_{j}x_{j}$).
\end{state}%

{\bf{}Proof:} We can assume $\varphi$ is not identically zero. First consider the possibility that
${\cal{}K}_{\varphi}=\emptyset$.  In this case we can write $\varphi$ as $\pm \varphi_{1}$, where
each component of $\varphi_{1}$ is a single-valued harmonic function on $\R^{n}$ which by the
relevant Liouville theorem is an affine function.

Thus we can assume that $\varphi$ is not identically zero and $0\in {\cal{}K}_{\varphi}$, and so
$$%
0<\sup_{\rho\in (0,\infty)}\rho^{1-n-2\alpha-2}\int_{\partial B_{\rho}}\varphi^{2}\le %
C\sup\rho^{-2-2\alpha}\max_{|x|=\rho}|\varphi|^{2} \le C[D\varphi]_{\alpha}^{2}<\infty.
$$%
Thus $\log(\rho^{1-n-2\alpha-2}\int_{\partial B_{\rho}}\varphi^{2})$ is bounded above on
$(0,\infty)$, and since it is a convex function of $t=\log\rho\in \R$ (by~2.3(1)), it must then be
constant. By Remark~2.4(1) this implies that $\varphi$ is homogeneous, and by~2.5 we conclude
that $\varphi$ is linear, contradicting $0\in {\cal{}K}_{\varphi}$. Thus ${\cal{}K}_{\varphi}\neq
\emptyset$ is impossible under the present hypotheses, and 2.6 is proved.

\section*{3\quad \bm{$C^{1,\alpha}$} estimates for a class of %
  linear equations.} %

Here we assume $\alpha\in (0,\fr 1/2)$ and $w\in C^{1,\alpha}(B_{1},\R^{k})$ is a 2-valued
symmetric function, and we continue to use the notation ${\cal{}K}_{w}=\{x\in B_{1}:w(x)=\{0,0\},
Dw(x)=\{0,0\}\}$.  Observe that then $|w(x)|\le C d(x)^{1+\alpha}$ and $|Dw(x)|\le C
d(x)^{\alpha}$ for some constant $C$, where $d(x)$ is distance of $x$ to ${\cal{}K}_{w}$.

We assume $B_{\rho}(y)\subset B_{1}$ and $w=(w^{1},\ldots,w^{k})$ satisfies a system of the form
$$%
\Delta w_{\kappa} + D_{j}(a^{ij}_{\kappa\lambda}D_{i}w_{\lambda}) + %
b_{\kappa}^{\lambda j}D_{j}w_{\lambda} %
+c^{\lambda}_{\kappa}w_{\lambda} = 0,\quad \kappa=1,\ldots,k %
\leqno{\bn{3.1}}
$$%
weakly on each ball $B_{\sigma}(z)\subset B_{\rho}(y)\setminus {\cal{}K}_{w}$, where the
coefficients $a^{ij}_{\kappa\lambda},b^{\lambda j}_{\kappa},c_{\kappa}^{\lambda}$ are single
valued. The main estimate is as follows:

\begin{state}{\bf{}3.2 Lemma.} %
  Let $\alpha\in(0,\delta(n))$ with $\delta(n)$ as in \emph{Lemma~2.5} and \emph{Corollary~2.6}. 
There is $\epsilon_{0}=\epsilon_{0}(\alpha,n,k)\in (0,1/2)$ such that the following holds for any
$\beta\ge 0$. Suppose~\emph{3.1} holds (weakly in $B_{\rho}(y)\setminus {\cal{}K}_{w}$) and
$$%
|a^{ij}_{\kappa\lambda}|_{0,B_{\rho}(y)} \le \epsilon_{0}, \quad %
\rho^{\alpha}[a^{ij}_{\kappa\lambda}]_{\alpha,B_{\rho}(y)} +\rho|b_{\kappa j}^{\lambda}|_{0} + %
\rho^{2}|c^{\lambda}_{\kappa}|_{0} \le \beta.
$$%
Then
 $$%
 \rho^{-1}\sup_{B_{\rho/2}(y)}|w|+
\rho^{-\alpha}\sup_{B_{\rho/2}(y)}|Dw|+[Dw]_{\alpha,B_{\rho/2}(y)} \le C
\rho^{-1-\alpha-n/2}\|w\|_{L^{2}(B_{\rho}(y))},
$$%
with $C=C(n,k,\beta)$.
\end{state}%

\medskip

{\textbf{Proof of Lemma~3.2:}} In view of scaling and standard interpolation inequalities it suffices
to consider the case $\rho=1$ and prove
 $$%
 [Dw]_{\alpha,B_{1/2}} \le C \|w\|_{L^{2}(B_{1})}. %
 \leqno{(1)}
$$%
To begin with we'll prove the formally weaker inequality that for each $\delta>0$ there is
$C=C(\delta,\beta,n)\ge 1$ and $\epsilon_{0}=\epsilon_{0}(\delta,\beta,n)>0$ such that the
hypotheses of the lemma imply
 $$%
 [Dw]_{\alpha,B_{1/2}} \le %
 \delta [Dw]_{\alpha,B_{1}}+ C (|w|_{0,B_{1}}+|Dw|_{0,B_{1}}). %
 \leqno{(2)}
$$%
If this fails then there are fixed $\delta,\beta>0$ and a sequence $w_{\ell}$ of solutions of~3.1, with
$a^{ij}_{\kappa\lambda,\ell},b^{\lambda j}_{\kappa,\ell},c^{\lambda}_{\kappa,\ell}$ in place of
$a^{ij}_{\kappa\lambda},b^{\lambda j}_{\kappa},c^{\lambda}_{\kappa}$ respectively, where
$$%
|a^{ij}_{\kappa\lambda,\ell}|_{0,B_{1}} \le \ell^{-1}, %
\quad [a^{ij}_{\kappa\lambda,\ell}]_{\alpha,B_{1}}+|b_{\kappa,\ell}^{\lambda j}|_{0,B_{1}} + %
|c_{\kappa,\ell}^{\lambda}|_{0,B_{1}} \le \beta, %
\leqno{(3)}
$$%
yet
$$%
[Dw_{\ell}]_{\alpha,B_{1/2}} > \delta[Dw_{\ell}]_{\alpha,B_{1}}+ %
C_{\ell} (|w_{\ell}|_{0,B_{1}} + |Dw_{\ell}|_{0,B_{1}}), \quad C_{\ell}\to \infty. %
\leqno{(4)}
$$%
Let ${\cal{}K}_{\ell}=\{x\in B_{1}:w_{\ell}(x)=\{0,0\}, Dw_{\ell}(x)=\{0,0\}\}$ and select distinct
points $x_{\ell},y_{\ell}\in B_{1/2}$ with
$|x_{\ell}-y_{\ell}|^{-\alpha}|Dw_{\ell}(x_{\ell})-Dw_{\ell}(y_{\ell})|>
\ha[Dw_{\ell}]_{\alpha,B_{1/2}}$. Observe that
$$%
\rho_{\ell}=|x_{\ell}-y_{\ell}|\to 0 %
\leqno{(5)}
$$%
because otherwise
$[Dw_{\ell}]_{\alpha,B_{1/2}}<2\rho_{\ell}^{-\alpha}|Dw_{\ell}(x_{\ell})-Dw_{\ell}(y_{\ell})|\le
4\sigma^{-\alpha}|Dw_{\ell}|_{0,B_{1}}$ with some fixed constant $\sigma>0$ (independent of
$\ell$), hence $[Dw_{\ell}]_{\alpha,B_{1/2}}<
4\sigma^{-\alpha}C_{\ell}^{-1}[Dw_{\ell}]_{\alpha,B_{1/2}}$ by~(4), a contradiction for $\ell$ large
enough to ensure $4\sigma^{-\alpha}C_{\ell}^{-1}<1$.

Then consider the possibilities:

\quad Case~1: \,\, $\Fr {\dist( \overline{x_{\ell}y_{\ell}},{\cal{}K}_{\ell})} /{|x_{\ell}-y_{\ell}|}$ is
bounded above

\quad Case~2: \,\, $\Fr {\dist( \overline{x_{\ell}y_{\ell}},{\cal{}K}_{\ell})} /{|x_{\ell}-y_{\ell}|}$ is not
bounded above,

where $\overline{x_{\ell}y_{\ell}}$ is the line segment joining $x_{\ell}$ and $y_{\ell}$.

In either case we define
$$%
\what w_{\ell}(x) = \sigma_{\ell} \rho_{\ell}^{-1-\alpha}w_{\ell}(y_{\ell}+\rho_{\ell}x), %
\quad x\in B_{\rho_{\ell}^{-1}/2}, %
$$%
where
$$%
\sigma_{\ell} = [Dw_{\ell}]_{\alpha,B_{1/2}}^{-1}
$$%
and observe that then $\what w_{\ell}$ is 2-valued symmetric $C^{1,\alpha}$ solution of a system of
the same form on $B_{\rho_{\ell}^{-1}/2}$ with
\begin{align*}%
  &\what w_{\ell}(0) = \rho_{\ell}^{-1-\alpha}\sigma_{\ell}w_{\ell}(y_{\ell}) , %
  \quad \what w_{\ell}(\xi_{\ell}) = \rho_{\ell}^{-1-\alpha}\sigma_{\ell}w_{\ell}(x_{\ell}) %
  \text{ where } \xi_{\ell}=\rho_{\ell}^{-1}(x_{\ell}-y_{\ell})\in S^{n-1},\\ %
  &[D\what w_{\ell}]_{\alpha,B_{\rho_{\ell}^{-1}/2}} <\delta^{-1}, \quad |D\what
w_{\ell}(\xi_{\ell})-D\what w_{\ell}(0)| \ge \ha. %
\end{align*}%
If Case~1 holds there is a bounded sequence $z_{\ell}$ with $z_{\ell}\in{\cal{}K}_{\what w_{\ell}}$
for each $\ell$. Thus we have $[D\what w_{\ell}]_{\alpha,B_{\rho_{\ell}^{-1}/2}}<\delta^{-1}$,
$|D\what w_{\ell}(\xi_{\ell})-D\what w_{\ell}(0)|\ge \ha$, $\what w_{\ell}(z_{\ell})=0$ and $D\what
w_{\ell}(z_{\ell})=0$, so in particular for each $x\in\R^{n}$ we have, for $\ell$ such that
$\rho_{\ell}^{-1}/2>|x|$, $|D\what w_{\ell}(x)|=|D\what w_{\ell}(x)-D\what w_{\ell}(z_{\ell})|\le
|x-z_{\ell}|^{\alpha}[D\what w_{\ell}]_{\alpha,B_{\rho_{\ell}^{-1}}}\le \delta^{-1}|x-z_{\ell}|^{\alpha}$
and $|\what w_{\ell}(x)|=|\what w_{\ell}(x)-\what w_{\ell}(z_{\ell})|\le
\delta^{-1}|x-z_{\ell}|^{1+\alpha}$.  By the Arzela-Ascoli theorem we can thus take a subsequence
such that $\what w_{\ell}$ converges locally in $C^{1}$ on $\R^{n}$ to a $C^{1,\alpha}$ symmetric
harmonic function $\varphi$ on $\R^{n}$ and which has the properties that
$[D\varphi]_{\alpha,\R^{n}}<\infty$ and $D\varphi$ is not constant.  By~2.6 we conclude that this
function must be affine if $\alpha\in (0,\delta(n))$, which contradicts the fact that $D\varphi$ is
non-constant.

In Case~2 we have a subsequence of $\ell$ and a corresponding sequence $R_{\ell}\to \infty$ with
$B_{R_{\ell}}\cap \{x\in \R^{n}:\what w_{\ell}(x)=0\text{ and }D\what w_{\ell}(x)=0\}=\emptyset$
and so $\what w_{\ell}$ gives a single-valued $C^{1,\alpha}$ functions $\wtilde w_{\ell}$ on
$B_{R_{\ell}}$ with $[D\wtilde w_{\ell}]_{\alpha,B_{R_{\ell}}}\le 1$.  Let $\wbar w=\wtilde
w_{\ell}-\wtilde w_{\ell}(0)-\sum_{j=1}^{n}x_{j}D_{j}\wtilde w_{\ell}(0)$ and $|D\wbar
w_{\ell}(\xi_{\ell})-D\wbar w_{\ell}(0)|\ge \ha$.  
Evidently $\wbar w_{\ell}$ weakly satisfies the system
$D_{j}((\delta_{ij}\delta_{\kappa\lambda}+a_{\kappa\lambda}^{ij}(y_{\ell}+\rho_{\ell}x))D_{i}\wbar
w_{\ell})=D_{j}((a_{\kappa\lambda}^{ij}(y_{\ell}+\rho_{\ell}x)-
a_{\kappa\lambda}^{ij}(y_{\ell}))D_{i}\wtilde w_{\ell}(0))=0$.  So since $|D\wtilde
w_{\ell}(0)|=\sigma_{\ell}\rho_{\ell}^{-\alpha}|Dw_{\ell}(y_{\ell})|\le
C_{\ell}^{-1}\rho_{\ell}^{-\alpha}$ by~(4), and $\rho_{\ell}^{-\alpha}
|a_{\kappa\lambda}^{ij}(y_{\ell}+\rho_{\ell}x)- a_{\kappa\lambda}^{ij}(y_{\ell})|\le
[a_{\kappa\lambda}^{ij}]_{\alpha}\le \beta$, we see the $\wbar w_{\ell}$ has a subsequence
converging locally in $C^{1}$ on $\R^{n}$ to a single-valued harmonic function $\varphi$ with
$\varphi(0)=0,D\varphi(0)=0$, $D\varphi$ non-constant and $[D\varphi]_{\alpha}\le 1$. But then,
since $\alpha<1$, $\varphi$ is linear, a contradiction.

Thus~(2) is proved, and by standard interpolation it implies
$$%
[Dw]_{\alpha,B_{1/2}} \le %
2 \delta [Dw]_{\alpha,B_{1}}+ C\|w\|_{L^{2}(B_{1})}, %
$$%
with $C=C(n,k,\delta,\beta)$, and the same must hold with appropriately scaled quantities over and
sub-ball of $B_{1}$; specifically,
\begin{align*}%
  \rho^{1+\alpha+n/2} [Dw]_{\alpha,B_{\rho/2}(y)} &\le %
  2 \delta \rho^{1+\alpha+n/2}[Dw]_{\alpha,B_{\rho}(y)}+ C\|w\|_{L^{2}(B_{\rho}(y))}\\ %
  &\le 2 \delta \rho^{1+\alpha+n/2}[Dw]_{\alpha,B_{\rho}(y)}+ C\|w\|_{L^{2}(B_{1})} %
\end{align*}%
for every ball $B_{\rho}(y)\subset B_{1}$, with $C=C(n,k,\delta,\beta)$, and, by taking
$\delta=\delta(n)$ suitably small and using a standard covering argument (see
e.g.~\cite[p398]{Sim97}), this implies
$$%
[Dw]_{\alpha,B_{1/2}} \le C\|w\|_{L^{2}(B_{1})} %
$$%
as claimed.

\section*{4\quad \bm{$\bn{2}$}-valued \bm{$C^{1}$} harmonic
functions---Part II}

Using 3.2 and the frequency function, we can now prove the following regularity theorem for
2-valued harmonic functions.

\begin{state}{\bf{}4.1 Lemma.} %
  Suppose $\varphi$, not identically zero, is a 2-valued symmetric $C^{1,\alpha}$ function on
$B_{1}$ which is smooth harmonic on $B_{1}\setminus {\cal{}K}_{\varphi}$.  Then the Hausdorff
dimension of ${\cal{}K}_{\varphi}$ is $\le n-2$, with equality if ${\cal{}B}_{\varphi}\neq \emptyset$,
and $\varphi$ is of class $C^{1,1/2}$, $|\varphi(x)|\le \smash{C\|\varphi\|_{L^{2}(B_{1})}d(x)^{3/2}}$
and $|D\varphi(x)|\le \smash{C\|\varphi\|_{L^{2}(B_{1})}d(x)^{1/2}}$ for all $x\in B_{1/2}$, where
$d(x)=\op{dist}(x,{\cal{}K}_{\varphi})$.

Furthermore, if $\varphi$ extends to all of\/ $\R^{n}$ as a homogeneous degree $\fr 3/2$ function
then, modulo composition with an orthogonal transformation of\/ $\R^{n}$, we have
$\varphi(x)=\varphi(re^{i\theta},0)=\{\pm Cr^{3/2}\sin 3\theta/2\}$ for some constant $C$, where
$r=(x_{1}^{2}+x_{2}^{2})^{1/2}$.
\end{state}%

{\textbf{Proof of Lemma~4.1:}} As in the proof of~2.5, we can assume that $k=1$, so that $\varphi$
is real-valued symmetric.

In view of the $C^{1,\alpha}$ estimates and Remark~2.3(4) we have
$$%
{\cal{}N}_{\varphi}(y)\ge 1+\alpha %
\leqno{(1)}
$$%
for each $y\in {\cal{}K}_{\varphi}$.

We are going to use a dimension reduction argument based on the monotonicity of the frequency
$N_{\varphi}(y,\rho)$ analogous to the procedure in \cite{Alm00}. To prove dimension
${\cal{}K}_{\varphi}\le n-2$ we suppose the contrary, that there is $s > n-2$ with ${\cal{}H}^{s} \,
({\cal{}K}_{\varphi} \cap B_{1}) > 0$. Let $\mu_{s}$ be the outer measure defined on subsets of $\R
^{n}$ by $\mu_{s}(A) = \inf \, \sum_{j=1}^{\infty} \rho_{j}^{s},$ where the $\inf$ is taken over
countable unions of balls $\cup_{j=1}^{\infty} \, B_{{\rho}_{j}}(y_{j})$ such that $A \subset
\cup_{j=1}^{\infty} \, B_{{\rho}_{j}}(y_{j})$. Note that ${\cal{}H}^{s} \, (A) > 0$ if and only if
$\mu_{s} \, (A) > 0$. Thus $\mu_{s} \, ({\cal{}K}_{\varphi} \cap B_{1}) > 0$, and hence there is a
point $z \in {\cal K}_{\varphi} \cap B_{1}$ such that the upper density $\theta_{\mu_{s}}^{\star} \,
({\cal K}_{\varphi}, z)=\limsup_{\rho\downarrow 0}\rho^{-s}\mu({\cal{}K}_{\varphi}\cap B_{\rho}(z))
>0$, so that there exists a sequence of positive numbers $\sigma_{j}\downarrow 0$ such that
$$%
\lim_{j \to \infty} \, \sigma_{j}^{-s}\mu_{s}\,%
({\cal{}K}_{\varphi} \cap B_{\sigma_{j}}(z)) >0 %
\leqno{(2)}
$$%
Set $\varphi_{j}(x) =(\sigma_{j}^{-n/2}\|\varphi\|_{L^{2}(B_{\sigma_{j}}(z))})^{-1}
\varphi(z+\sigma_{j}x)$ for $x\in B_{\sigma_{j}^{-1}(1-|z|)}$.  In view of~(2) we have
$$%
\liminf_{j\to\infty }\mu_{s}({\cal{}K}_{\varphi_{j}}\cap \overline B_{1}) >0 \leqno{(3)}
$$%
By virtue of the $C^{1,\alpha}$ estimates of~3.2, and the Remarks~2.3(4),(5) and~2.4(1) we also
have that a subsequence $\varphi_{j^{\prime}}$ converges locally in $C^{1}$ on $\R^{n}$ to a
2-valued symmetric $C^{1,\alpha}\cap W^{2,2}$ function $\psi$ which is homogeneous of some
degree $>1$ (in fact $\ge 1+\alpha$) on $\R^{n}$, $\|\psi\|_{L^{2}(B_{1})}=1$, harmonic on
$B_{1}\setminus{\cal{}K}_{\psi}$, $0\in{\cal{}K}_{\psi}$ and
${\cal{}N}_{\psi}(0)={\cal{}N}_{\varphi}(z)$.  Also $\mu_{s}({\cal{}K}_{\psi}\cap\ovl{B}_{1}) \ge
\limsup_{j^{\prime}\to\infty}\mu_{s}({\cal{}K}_{\varphi_{j^{\prime}}}\cap \ovl{B}_{1})(>0$ by~(3)),
because if $\delta>0$ we can choose $B_{\rho_{\ell}}(y_{\ell})$ with
${\cal{}K}_{\psi}\cap\ovl{B}_{1}\subset\cup_{\ell}B_{\rho_{\ell}}(y_{\ell})$ and
$\sum_{\ell}\rho_{\ell}^{s}<\mu_{s}({\cal{}K}_{\psi}\cap\ovl{B}_{1})+\delta$, and, since the
convergence is $C^{1}$, we see that ${\cal{}K}_{\varphi_{j^{\prime}}}\cap \ovl{B}_{1}\subset
\cup_{\ell}B_{\rho_{\ell}}(y_{\ell})$ for all sufficiently large $j^{\prime}$ and hence
$\mu_{s}({\cal{}K}_{\psi}\cap\ovl{B}_{1})>\mu_{s}({\cal{}K}_{\varphi_{j^{\prime}}}\cap \ovl
B_{1})-\delta$ for all sufficiently large $j^{\prime}$, whence $\mu_{s}({\cal{}K}_{\psi}\cap\ovl B_{1})
\ge \limsup_{j^{\prime}\to\infty}\mu_{s}({\cal{}K}_{\varphi_{j^{\prime}}}\cap \ovl B_{1})>0$ as
claimed.

Thus we can repeat the argument with $\psi$ in place of $\varphi$ and with new base point $z\in
{\cal{}K}_{\psi}\setminus\{0\}$ as origin, and this gives us a homogeneous harmonic $\psi_{1}$, not
identically zero, but still with $\mu_{s}({\cal{}K}_{\psi_{1}})>0$ and with $0\in{\cal{}K}_{\psi_{1}}$
and $\psi_{1}$ invariant under composition with translations in the direction of $\lambda z$
($\lambda\in \R$).  Since $s>n-2$ we can repeat this process a further $n-2$ times to give finally a
homogeneous $C^{1,\alpha}$ harmonic $\psi_{n-1}$, not identically zero and such that there is an
$(n-1)$-dimensional subspace $L\subset{\cal{}K}_{\psi_{n-1}}$ with $\psi_{n-1}$ invariant under
composition with translations in the direction of vectors in $L$. However then, assuming without
loss of generality that $L=\R^{n-1}\times\{0\}$, we have $\psi_{n-1}(x_{1},\ldots,x_{n})=f(x_{n})$
(independent of $x_{1},\ldots,x_{n-1}$) with $f$ a 2-valued homogeneous symmetric $C^{1}$
function, harmonic on $\R\setminus\{0\}$, not identically zero, and $f(0)=f^{\prime}(0)=\{0,0\}$. 
Clearly no such function $f$ exists so we have a contradiction. Hence ${\cal{}K}_{\varphi}$ is of
Hausdorff dimension $\le n-2$ as claimed.

Now suppose $y\in {\cal{}B}_{\varphi}$. Then ${\cal{}K}_{\varphi}\cap B_{\rho}(y)$ has positive
$(n-2)$-dimensional Hausdorff measure for each $\rho\in (0,1-|y|)$, because otherwise (see the
appendix) $B_{\rho}(y)\setminus{\cal{}K}_{\varphi}$ would be simply connected for some $\rho\in
(0,1-|y|)$ and we could write $\varphi|B_{\rho}(y)=\{\pm\varphi_{1}\}$, where $\varphi_{1}$ is a
single-valued harmonic function on $B_{\rho}(y)$, contradicting the fact that
$y\in{\cal{}B}_{\varphi}$.  Thus we can repeat the above dimension reducing argument with
$s=n-2$ and starting at points $z$ arbitrarily close to $y$, obtaining, after a total of $n-2$ steps and
modulo a composition with an orthogonal transformation of $\R^{n}$, a homogeneous 2-valued
symmetric $C^{1,\alpha}$ harmonic function $\psi$, with $\|\psi\|_{L^{2}(B_{1})}=1$,
$0\in{\cal{}K}_{\psi}$, ${\cal{}N}_{\varphi}(z)\ge {\cal{}N}_{\psi}(0)$, $\psi$ homogeneous of degree
$m>1$ and invariant under composition with any translation in the direction of any vector in
$\{0\}\times\R^{n-2}$, so that $\psi$ is a function of just the first 2-variables
$x_{1}=r\cos\theta,x_{2}=r\sin\theta$.  Then (modulo composition with a rotation of $\R^{2}$) we
must be able to write $\psi=\{\pm Cr^{m/2}\cos m\theta/2\}$ with $m$ an integer $\ge 3$, and
hence ${\cal{}N}_{\!  \varphi}(z)\ge {\cal{}N}_{\!  \psi}(0)\ge 3/2$.  Since we can repeat this argument
at a sequence of points $z\in{\cal{}K}_{\varphi}$ with $z\to y$, we conclude by the upper
semicontinuity~2.3(3) of ${\cal{}N}_{\varphi}$ that ${\cal{}N}_{\!  \varphi}(y)\ge 3/2$.  Of course at
points $y\in {\cal{}K}_{\varphi}\setminus {\cal{}B}_{\varphi}$ we must have ${\cal{}N}_{\! 
\varphi}(y)\ge 2> 3/2$.  Thus in any case ${\cal{}N}_{\!  \varphi}(y)\ge 3/2$ for
$y\in{\cal{}K}_{\varphi}$ and hence by Remark~2.3(4) and standard estimates for single-valued
harmonic functions (in balls $B_{\rho}(z)$ with $B_{\rho}(z)\cap {\cal{}K}_{\varphi}=\emptyset$) we
have $|\varphi|\le C \|\varphi\|_{L^{2}(B_{1})}d^{3/2}$ and $|D\varphi|\le
C\|\varphi\|_{L^{2}(B_{1})}d^{1/2}$ on $B_{1/2}$, as claimed.

Of course if $\varphi$ is homogeneous of degree $3/2$ on $\R^{n}$ to begin with, then
$0\in{\cal{}B}_{\varphi}$ (otherwise $\varphi$ would be $\pm$ a single valued harmonic function,
which of course would have integer homogeneity), the first step of the above dimension reducing
argument would not be needed, and, by Remarks~2.4(1),(2), $\varphi$ itself would be invariant
under composition with translations by vectors in some $(n-2)$-dimensional subspace; that is,
modulo composition with a rotation of coordinates, we would have
$\varphi(x_{1},x_{2},y)=\varphi(re^{i\theta},0)=\{\pm Cr^{3/2}\sin 3\theta/2\}$ for some constant
$C$, as claimed.

\medskip

We next establish a gap lemma analogous to~2.5 for orders of homogeneity $\in (\fr 3/2,\fr 3/2 +
\delta)$.

\begin{state}{\bf{}4.2 Lemma.} %
  There is $\delta=\delta(n)\in (0,1)$ such that if $\varphi(x)=|x|^{\sigma}\varphi(|x|^{-1}x)$ is
2-valued symmetric locally $C^{1,\alpha}$ homogeneous degree $\sigma$ function with $\sigma\in
(\fr 3/2,\fr 3/2 +\delta)$ and with $\varphi$ harmonic on $\R^{n}\setminus {\cal{}K}_{\varphi}$,
then $\varphi\equiv \{0,0\}$ on $\R^{n}$.
\end{state}%

{\textbf{Proof:}} Again, as in the proof of~2.5, we can assume that $k=1$, so that $\varphi$ is
real-valued symmetric.  The proof in this case is similar to Lemma~2.5, but not quite as simple. If
the result is false, we would have a sequence $\varphi_{\ell}$ of 2-valued symmetric $C^{1,\alpha}$
functions $\varphi_{\ell}$ harmonic on $\R^{n}\setminus{\cal{}K}_{\varphi_{\ell}}$, $\varphi_{\ell}$
not identically zero, and homogeneous of degree $\sigma_{\ell}$, where $\fr 3/2<
\sigma_{\ell}\downarrow \fr 3/2$. Assume that $\varphi_{\ell}$ is normalized so that
$\|\varphi_{\ell}\|_{L^{2}(B_{1})}=1$. Then, by~3.2, $\varphi_{\ell}\to \varphi$ in $C^{1}$, where
$\varphi$ is $C^{1,\alpha}$ and homogeneous of degree 3/2, and, by~4.1, modulo an orthogonal
change of independent variables and a rescaling we have $\varphi(x)= \varphi(re^{i\theta},0)=\{\pm
r^{3/2}\sin3\theta/2\}$, assuming $x=(re^{i\theta},x_{3},\ldots,x_{n})\in \R^{2}\times \R^{n-2}$ and
$r=(x_{1}^{2}+x_{2}^{2})^{1/2}$.  Take ${\cal{}K}_{\ell}=\{x\in\R^{n}:
|\varphi_{\ell}(x)|=0,|D\varphi_{\ell}(x)|=0\}$.

Next let $\varphi_{\ell,j}=D_{j}\varphi_{\ell}$, $j=1,2,3$, on $\R^{n}\setminus {\cal{}K}_{\ell}$. 
Direct computation shows that $\varphi_{\ell,1}\to \{\pm c r^{1/2}\sin\theta/2\}$ and
$\varphi_{\ell,2}\to -\{\pm c r^{1/2}\cos\theta/2\}$ (for suitable $c=c(n)>0$).

Now choose $c_{\ell},d_{\ell}$ to ensure $\varphi_{\ell,3}
-c_{\ell}\varphi_{\ell,1}-d_{\ell}\varphi_{\ell,2}$ is orthogonal to
$\spn\{\varphi_{\ell,1},\varphi_{\ell,2}\}$ on $S^{n-1}\cap K$---i.e.\ (Cf.\ (4) in the proof of
Lemma~2.5)
$$%
\sint_{S^{n-1}\cap K}(\varphi_{\ell,3} %
-c_{\ell}\varphi_{\ell,1}-d_{\ell}\varphi_{\ell,2})\varphi_{\ell,j}=0,\quad j=1,2, %
\leqno{(1)}
$$%
where $K=\{(x_{1},x_{2},y)\in\R^{2}\times\R^{n-2}:r^{2}\le \fr 1/2 |y|^{2}\}$, and notice that we can
assume that $\varphi_{\ell,3}-c_{\ell}\varphi_{\ell,1}-d_{\ell}\varphi_{\ell,2}$ are not identically zero
on $\R^{n}$---otherwise the directional derivative of $\varphi_{\ell}$ in the direction
$\eta_{\ell}=(-c_{\ell},-d_{\ell},1,0)$ would be zero on $\R^{n}$, and this would imply that
$\varphi_{\ell}$ is cylindrical in the direction $\eta_{\ell}$, which for $n\ge 2$ would enable us to
reduce the dimension $n$ in the statement of the lemma (from $n$ to $n-1$); and so in fact, since
the only homogeneous harmonic functions on $\R$ are the linear functions which are homogeneous
of degree~1, we can assume without loss of generality that this does not happen.  So we can define
$\psi_{\ell}=\|\varphi_{\ell,3}- c_{\ell}\varphi_{\ell,1}-d_{\ell}\varphi_{\ell,2}\|_{L^{2}( S^{n-1})}^{-1}
(\varphi_{\ell,3} -c_{\ell}\varphi_{\ell,1}-d_{\ell}\varphi_{\ell,2})$.

For $\delta>0$ let $\gamma_{\delta}(t)=\op{sgn} (t)\,\max\{|t|-\delta,0\}$ and observe that then
$\gamma_{\delta}(\psi_{\ell})$ has compact support in $S^{n-1}\setminus {\cal{}K}_{\ell}$, and
hence (since $\psi_{\ell}$ is harmonic and homogeneous degree $\sigma_{\ell}-1$) we have
$\int_{S^{n-1}}\nabla_{S^{n-1}}\psi_{\ell}\cdot \nabla_{S^{n-1}}\gamma_{\delta}(\psi_{\ell})=
(\sigma_{\ell}-1)(\sigma_{\ell}+n-3)\int_{ S^{n-1}} \psi_{\ell}\gamma_{\delta}(\psi_{\ell})$.  Thus
letting $\delta\downarrow 0$ we obtain the identity $\int_{S^{n-1}}|\nabla_{ S^{n-1}}
\psi_{\ell}|^{2}=(\sigma_{\ell}-1)(\sigma_{\ell}+n-3)\int_{S^{n-1}}\psi_{\ell}^{2}=
(\sigma_{\ell}-1)(\sigma_{\ell}+n-3)$ for each $\ell=1,2,\ldots$, and in particular we have a uniform
bound on the $L^{2}$-norm of $\nabla_{ S^{n-1}} \psi_{\ell}$ and so by Rellich's theorem in the
limit as $\ell\to \infty$ this gives us 2-valued symmetric homogeneous degree $\ha$, harmonic $\psi$
in $W^{1,2}_{\text{loc}}(\R^{n}\setminus (\{0\}\times\R^{n-2}))$, and (by lower semicontinuity of
the norm with respect to weak $L^{2}$ convergence)
$$%
\sint_{S^{n-1}}|\nabla \psi|^{2} \le %
\fr 1/2(n-\fr 3/2)\sint_{S^{n-1}}\psi^{2}=\fr 1/2(n-\fr 3/2). %
\leqno{(2)}
$$%
Since $\psi(r,\theta,y)$ is locally the $C^{0}$ limit of $\psi_{\ell}(r,\theta,y)$, each of which can be
viewed for each fixed $r,y$ with $r>0$ and $r^{2}+|y|^{2}=1$ as a $4\pi$-periodic single valued
function of $\theta$, so $\psi(r,\theta,y)$ can also be viewed as a smooth $4\pi$-periodic single
valued function of $\theta$ (with $\psi(r,\theta+2\pi,y)=-\psi(r,\theta,y)$ by the symmetry), and we
see that~(1) above ensures that
$$%
\psi \text{ is orthogonal to } r^{1/2}\cos\theta/2,\,r^{1/2}\sin \theta/2 \,\, %
(\theta\in [0,4\pi]) \text{ in }L^{2}(S^{n-1}\cap K) %
\leqno{(3)}
$$%
(and the statement~(3) makes sense).  Defining $\psi_{0}=r^{-1/2}\psi$, then $\psi_{0}$ is 2-valued
symmetric homogeneous of degree zero on $\R^{n}\setminus (\{0\}\times\R^{n-2})$ and, with
$\nabla = $ gradient on $S^{n-1}$,
\begin{align*}%
  \sint_{S^{n-1}}|\nabla \psi|^{2} %
  &= \sint_{S^{n-1}}|\nabla (r^{1/2}\psi_{0})|^{2} \\ %
  &= \sint_{S^{n-1}}\Bl(r|\nabla\psi_{0}|^{2}+ %
  \fr 1/4 r^{-1}|\nabla r|^{2}\psi_{0}^{2}+2r^{1/2}\psi_{0}\nabla(r^{1/2})\cdot %
  \nabla \psi_{0}\Br) \\ %
  &= \sint_{S^{n-1}}\Bl(r|\nabla \psi_{0}|^{2}+\fr 1/4 r^{-1}|\nabla %
  r|^{2}\psi_{0}^{2}+\fr 1/2 \nabla(r)\cdot\nabla (\psi_{0}^{2})\Br) \\ %
  &= \sint_{S^{n-1}}\Bl(r|\nabla\psi_{0}|^{2}+\fr 1/4 r^{-1}|\nabla r|^{2}\psi_{0}^{2}- %
  \fr 1/2 \psi_{0}^{2}\Delta_{S^{n-1}}r\Br), %
\end{align*}%
and since $\Delta_{S^{n-1}}r={\fr 1/r}-(n-1)r$ and $|\nabla r|^{2}=1-r^{2}$, this implies
$$%
\sint_{S^{n-1}}|\nabla \psi|^{2} = \sint_{S^{n-1}}{\Bl(r|\nabla \psi_{0}|^{2}+ %
  \Bl(\fr 1/2\Bl(n-{\fr 3/2}\Br)r -{\fr 1/4}r^{-1}\Br)\psi_{0}^{2}\Br)}.  %
\leqno{(4)}
$$%
Since $\psi_{0}$ is homogeneous of degree zero we have, writing points in $\R^{n}$ as
$(x,y)=(re^{i\theta},y)$ with $r>0$ and $y\in \R^{n-2}$, $|\nabla \psi_{0}|^{2}=|D\psi_{0}|^{2}=
r^{-2}({\fr
  \partial/{\partial\theta}}\psi_{0})^{2}+({\fr \partial/{\partial r}}\psi_{0})^{2}+|D_{y}\psi_{0}|^{2}$ on
$\R^{n}\setminus(\{0\}\times\R^{n-2})$, and so
$$%
\sint_{S^{n-1}}r|\nabla \psi_{0}|^{2}\ge %
\sint_{S^{n-1}}r^{-1}\Bigl({\Fr {\partial\psi_{0}}/{\partial\theta}}\Bigr)^{2} %
\leqno{(5)}
$$%
with equality if and only if $\psi_{0}$ is cylindrical in the sense that $\psi_{0}$ is independent of the
variables $r,y$---i.e., equality holds in~(5) if and only if $\psi_{0}$ is a function of $\theta$ alone. 
On the other hand we have the general fact (using Fourier series for $f,f'$) that
$\int_{0}^{4\pi}(f^{\prime}(\theta))^{2}\,d\theta \ge {\fr 1/4}\int_{0}^{4\pi}f^{2}(\theta)\,d\theta$ for
any Lipschitz $4\pi$-periodic function $f$ with $f(\theta+2\pi)\equiv -f(\theta)$ (so that $f$ defines
a 2-valued symmetric function on $S^{1}$), and equality holds if and only if $f(\theta)$ has the
form $a\cos\theta/2+b\sin\theta/2$ for some constants $a,b$, and so we deduce from~(5) that
$$%
\sint_{S^{n-1}}r|\nabla \psi_{0}|^{2} \ge {\fr 1/4}\sint_{S^{n-1}}r^{-1}\psi_{0}^{2}, %
\leqno{(6)}
$$%
with equality if and only if $\psi_{0}\equiv \{\pm(a\cos\theta/2+b\sin\theta/2)\}$ for some constants
$a,b$. Thus finally by~(4),(6)
$$%
\sint_{S^{n-1}}|\nabla\psi|^{2}\ge %
\sint_{S^{n-1}}\fr 1/2 \Bl(n-{\fr 3/2}\Br)r\psi_{0}^{2} = %
\sint_{S^{n-1}}\fr 1/2 \Bl(n-{\fr 3/2}\Br)\psi^{2}=\fr 1/2 \Bl(n-{\fr 3/2}\Br), %
$$%
with equality if and only if $\psi\equiv \pm r^{1/2}(a\cos\theta/2+b\sin\theta/2)$ for some constants
$a,b$. However by~(2) we do have equality, and hence $\psi(r,\theta,y)$ has the form
$r^{1/2}(a\cos\theta/2+b\sin\theta/2)$ ($\theta\in [0,4\pi]$) for some constants $a,b$, which
contradicts~(3).

\section*{5\quad Regularity for \bm{$u_{a}= %
    \fr 1/2 (u_{1}+u_{2})$} and %
  \bm{$v=\pm\ha(u_{1}-u_{2})$}, Part~I}

First we observe that the 2-valued symmetric $C^{1,\alpha}(B_{1})$ function
$v=\{\pm\ha(u_{1}-u_{2})\}$ is just $(u-\{u_{a},u_{a}\})$ (so we can write $u=\{u_{1},u_{2}\}=u_{a}+
v=\{u_{a}\pm \ha (u_{1}-u_{2})\}$) and $v$ and $Dv$ vanish to orders $1+\alpha$ and $\alpha$
respectively on ${\cal{}K}_{u}=\{x:u_{1}(x)=u_{2}(x)\text{ and }Du_{1}(x)=Du_{2}(x)\}$; thus for
each $\theta\in (0,1)$
$$%
|v(x)|\le Cd(x)^{1+\alpha},\quad |Dv(x)|\le Cd(x)^{\alpha},\quad x\in B_{\theta} %
\leqno{\hbox{\bf 5.1}}
$$%
for some constant $C=C(\theta)$, where, here and subsequently, $d(x)=\dist(x,{\cal{}K}_{v})$.

Before we begin the proof of the main estimates we need some general remarks about the nature of
the equations governing $u=\{u_{1},u_{2}\}$.

First, near points of $B_{1}\setminus {\cal{}K}_{u}$ we can represent $u$ as an ordered pair
$(u_{1},u_{2})$ with each $u_{\ell}=(u_{\ell}^{1},\ldots,u_{\ell}^{k})$ a smooth $\R^{k}$-valued
solution of the minimal surface system. Thus
$$%
\ssum_{i,j=1}^{n}D_{i}(G^{ij}(Du_{\ell})D_{j}u_{\ell}^{\kappa})=0, %
\quad \kappa=1,\ldots,k,\,\, \ell=1,2, %
\leqno{\hbox{\bf 5.2}}%
$$%
on $B_{1}\setminus{\cal{}K}_{u}$, where
\begin{align*}%
  &G^{ij}(p)=\sqrt{g(p)}g^{ij}(p), \text{ with $(g^{ij}(p))=(g_{ij}(p))^{-1}$, } \\ %
  &\hskip0.9in\text{$g_{ij}(p)=\delta_{ij}+ %
    \textstyle\sum_{\kappa=1}^{k}p_{i}^{\kappa}p_{j}^{\kappa}$ %
    and $g(p)=\det (g_{ij}(p))$, $p=(p_{i}^{\kappa})\in \R^{n}\otimes\R^{k}$.}  %
\end{align*}%
We also recall that the fact that the mean curvature of $\graph u$ is zero in $(B_{1}\setminus
{\cal{}K}_{u})\times \R^{k}$ implies the identities
$$%
\ssum_{i=1}^{n}D_{i}\bigl(G^{ij}(Du_{\ell})\bigr)=0,\quad j=1,\ldots,n, \ell=1,2%
\leqno{\hbox{\bf 5.3}}%
$$%
on $B_{1}\setminus{\cal{}K}_{u}$, so~5.2 can also be written in non-divergence form
$$%
\ssum_{i,j=1}^{n} G^{ij}(Du_{\ell})D_{i}D_{j}u_{\ell}^{\kappa} =0, %
\quad \kappa=1,\ldots,k,\,\, \ell=1,2, %
\leqno{\hbox{\bf 5.4}}%
$$%
on $B_{1}\setminus{\cal{}K}_{u}$.

We want to write these equations on $B_{1}\setminus {\cal{}K}_{u}$ in more readily usable form, in
terms of
$$%
u_{a} =\ha(u_{1}+u_{2}),\quad v=\{\pm\ha (u_{1}-u_{2})\}
$$%
(note that $u_{a}$ is then single-valued and $v$ is 2-valued symmetric) and the functions
$$%
\begin{aligned}%
  A^{ij}(p,q) &= G^{ij}(p+q) + G^{ij}(p-q) \\ E^{ij\ell}_{\lambda}(p,q) &=
\int_{-1}^{1}G^{ij}_{p^{\lambda}_{\ell}}(p+s q)\,ds .
\end{aligned}%
\leqno{\hbox{\bf 5.5}}
$$%
Observe that $A^{ij},E^{ij\ell}_{\lambda}$ are real-analytic functions of
$p,q\in\R^{n}\otimes\R^{k}$ and (using the definitions and the symmetry $G^{ij}(p)=G^{ij}(-p)$)
we have the symmetries
$$%
\begin{aligned}%
  A^{ij}(p,q) %
  &= A^{ij}(p,-q) = A^{ij}(-p,q) = A^{ij}(-p,-q) \\ %
  E^{ij\ell}_{\lambda}(p,q) &= E^{ij\ell}_{\lambda}(p,-q) = -E^{ij\ell}_{\lambda}(-p,q), %
\end{aligned}%
\leqno{\hbox{\bf 5.6}}
$$%
and in particular
$$%
E^{ij\ell}_{\lambda}(0,q) =0 \text{ and }D_{q}A^{ij}(p,0)=0, %
\leqno{\hbox{\bf 5.7}}
$$%
whence we can write
$$%
E^{ij\ell}_{\lambda}(p,q) =\sum_{\kappa=1}^{k}\sum_{h=1}^{n} %
E_{\kappa\lambda}^{ijh\ell}(p,q)p_{h}^{\kappa} %
$$%
for suitable real-analytic $E^{ijh\ell}_{\kappa\lambda}$ with the symmetries
$$%
E_{\kappa\lambda}^{ijh\ell}(p,q)= E_{\kappa\lambda}^{ijh\ell}(p,-q)= %
E_{\kappa\lambda}^{ijh\ell}(-p,-q).  %
\leqno{\hbox{\bf 5.8}}
$$%
Notice that $E^{ij\ell}_{\lambda}$ arises naturally from the calculus identity
$$%
G^{ij}(p+q) - G^{ij}(p-q) =\int_{-1}^{1}{d\over{}dt}G^{ij}(p+tq)\,dt= %
\sum_{\lambda=1}^{k}\sum_{\ell=1}^{n}E^{ij\ell}_{\lambda}(p,q)q^{\lambda}_{\ell}, %
\quad p,q\in\R^{n}\otimes\R^{k} %
$$%
In view of the evenness of $A^{ij}(p,q),E^{ij\ell}_{\lambda}(p,q),E^{ijh\ell}_{\kappa\lambda}(p,q)$
with respect to $q$ we see that
$$%
A^{ij}(Du_{a},Dv),\,E^{ij\ell}_{\lambda}(Du_{a},Dv),\,\, E^{ijh\ell}_{\kappa\lambda}(Du_{a},Dv)
$$%
are actually single-valued (rather than 2-valued) and~5.3 implies
$$%
\left\{%
\begin{aligned}%
  D_{i}(A^{ij}(Du_{a},Dv)) &=0 \\ %
  D_{i}(E^{ij\ell}_{\lambda}(Du_{a},Dv)D_{\ell}v^{\lambda})&=\{0,0\} %
\end{aligned}%
\right. %
\leqno{\hbox{\bf 5.3'}}
$$%
on $B_{1}\setminus{\cal{}K}_{u}$.  Also by first taking the difference of the two equations in~5.2
and writing $u_{1}=u_{a}+\ha (u_{1}-u_{2})$ and $u_{2}=u_{a}-\ha (u_{1}-u_{2})$ we obtain
\begin{align*}%
  \tag*{\bf 5.9} %
  &D_{i}\Bigl(\Bigl(G^{ij}\bigl(Du_{a}+ Dv\bigr)+G^{ij}\bigl(Du_{a}- Dv\bigr)\Bigr) %
  D_{j}v^{\kappa}\Bigr)\\ %
  &\hskip1.2in +D_{i}\Bigl(\Bigl(G^{ij} %
  \bigl(Du_{a}+Dv\bigr)-G^{ij}\bigl(Du_{a}- Dv\bigr)\Bigr)D_{j}u_{a}^{\kappa}\Bigr) %
  =\{0,0\} %
\end{align*}%
and, by taking the sum,
\begin{align*}%
  \tag*{\bf 5.9'} %
  &D_{i}\Bigl(\Bigl(G^{ij}\bigl(Du_{a}+ Dv\bigr)+G^{ij}\bigl(Du_{a}- Dv\bigr)\Bigr) %
  D_{j}u^{\kappa}_{a}\Bigr)\\ %
  &\hskip1.2in + D_{i}\Bigl(\Bigl(G^{ij} %
  \bigl(Du_{a}+ Dv\bigr)-G^{ij}\bigl(Du_{a}- Dv\bigr)\Bigr)D_{j}v^{\kappa}\Bigr) =0 %
\end{align*}%
on $B_{1}\setminus{\cal{}K}_{u}$.  In terms of the single-valued functions
$A^{ij}(Du_{a},Dv),E^{ij\ell}_{\lambda}(Du_{a},Dv)$, these equations can be written
\begin{align*}%
  \tag*{\bf 5.10}&D_{i}\Bigl(A^{ij}(Du_{a},Dv)D_{j}v^{\kappa} + %
  E^{ij\ell}_{\lambda}(Du_{a},Dv)D_{\ell}v^{\lambda}D_{j}u_{a}^{\kappa}\Bigr) %
  = \{0,0\}\\ %
  \tag*{\bf 5.11}& D_{i}\Bigl(A^{ij}(Du_{a},Dv)D_{j}u_{a}^{\kappa} + %
  E^{ij\ell}_{\lambda}(Du_{a},Dv)D_{\ell}v^{\lambda}D_{j}v^{\kappa}\Bigr)= 0 %
\end{align*}%
on $B_{1}\setminus {\cal{}K}_{u}$ for $\kappa=1,\ldots,k$.  In view of~5.3$'$ these can also be
written in non-divergence form as
\begin{align*}%
  \tag*{\bf 5.10'}&A^{ij}(Du_{a},Dv)D_{i}D_{j}v^{\kappa} + %
  E^{ij\ell}_{\lambda}(Du_{a},Dv)D_{\ell}v^{\lambda}D_{i}D_{j}u_{a}^{\kappa} %
  = \{0,0\}\\ %
  \tag*{\bf 5.11'}& A^{ij}(Du_{a},Dv)D_{i}D_{j}u_{a}^{\kappa} + %
  E^{ij\ell}_{\lambda}(Du_{a},Dv)D_{\ell}v^{\lambda}D_{i}D_{j}v^{\kappa}= 0 %
\end{align*}%

\smallskip

{\bf{}5.12 Remark:} We shall use the fact that~5.11 actually holds in the weak sense on all of
$B_{1}$ (i.e.\ across ${\cal{}K}_{u}$ also.) To see this, note that the first variation formula~1.2 for
the stationary varifold $G=\graph u$ can be written $\int_{G}\nabla_{i}\zeta\,d{\cal{}H}^{n}=0$ for
any $i=1,\ldots,n+k$ and any smooth $\zeta$ with compact support in $B_{1}\times\R^{k}$, where
$\nabla\zeta=(\nabla_{1}\zeta,\ldots,\nabla_{n+k}\zeta)=P_{x}(D\zeta)$, where $P_{x}$ is orthogonal
projection onto the tangent space of $G$ at $(x,u(x))$.  Since $G$ is bounded we can of course take
$\zeta$ to be a $C^{\infty}_{c}(B_{1})$ function (i.e.\ $\zeta$ can be taken independent of the
variables $x_{n+1},\ldots,x_{n+k}$), in which case for $i=1,\ldots,n$ the above formula can be written
$$%
\sint_{B_{1}}\ssum_{i=1}^{n}(\nu_{i}^{\kappa}(Du_{1})+\nu_{i}^{\kappa}(Du_{2})) D_{i}\zeta=0, %
\quad \zeta\in C^{1}_{c}(B_{1}), %
$$%
where $\nu_{i}^{\kappa}(Du_{\ell})=\sum_{j=1}^{n}G^{ij}(Du_{\ell})D_{j}u_{\ell}^{\kappa}$ as
in~5.2, so indeed this identity is exactly the weak form, over all of $B_{1}$, of 5.11.

\medskip

With these facts at our disposal we can prove some initial rough $C^{1, 1-\epsilon}$ estimates for
$u_{a}$ and $C^{1, \fr 1/2-\epsilon}$ estimates for $v$:

\begin{state}{\bf{}5.13 Lemma.} %
  There is $\epsilon_{0}(n)$ such that if~1.6 holds with $\epsilon_{0}\le \epsilon_{0}(n)$ and if
$u_{a}=\fr 1/2 (u_{1}+u_{2})$ and $v=\{\pm\ha (u_{1}-u_{2})\}$, then for each $\epsilon>0$ there is
$C=C(n,\epsilon)$ with
\al{%
  &\hskip1.2in %
  \sup_{x\in B_{1/2}\setminus {\cal{}K}_{u}} d(x)^{-\fr 1/2+\epsilon}|Dv(x)| <C\epsilon_{0}, \\ %
  &\sup_{y\in B_{1/2}\cap {\cal{}K}_{u},x\in B_{1/2}}|x-y|^{- 1+\epsilon} %
  |Du_{a}(x)-Du_{a}(y)| <C\epsilon_{0},\, %
  \sup_{x\in B_{1/2}}d(x)^{\epsilon} |D^{2}u_{a}(x)| <C\epsilon_{0}. %
}%
\end{state}%

{\bf{}Remark:} We ultimately show that the above inequalities are valid with $\epsilon=0$ but we
need the above lemma in the course of the proof of this.

\bigskip

{\bf{}Proof of Lemma~5.13:} We can assume $B_{1/2}\cap {\cal{}K}_{u}\neq \emptyset$, otherwise
we use standard quasilinear elliptic theory (for single-valued solutions) to give the claimed estimates
over $B_{1/2}$. We claim that for each $\epsilon>0$ there is $\rho_{0}=\rho_{0}(\epsilon)\in (0,\fr
1/4)$ such that $y\in B_{1/2}\cap {\cal{}K}_{u}$ and $\rho<\rho_{0}\Rightarrow
(\rho/2)^{-n/2}\|v\|_{L^{2}(B_{\rho/2}(y))}\le 2^{-3/2+\epsilon}\rho^{-n/2}\|v\|_{L^{2}(B_{\rho}(y))}$.
 Otherwise there would be sequences $y_{\ell}\in \ovl{B}_{1/2}\cap {\cal{}K}_{u}$ and
$\rho_{\ell}\downarrow 0$ such that
$$%
(\rho_{\ell}/2)^{-n/2}\|v\|_{L^{2}(B_{\rho_{\ell}/2}(y_{\ell}))}>
2^{-3/2+\epsilon}\rho_{\ell}^{-n/2}\|v\|_{L^{2}(B_{\rho_{\ell}}(y_{\ell}))}
$$%
for all $\ell$.  After a translation and rescaling transforming $B_{\rho_{\ell}}(y_{\ell})$ to
$B_{1}(=B_{1}(0))$ and $v(x)$ to $v_{\ell}(x)=\beta_{\ell}^{-1}v(y_{\ell}+\rho_{\ell}x)$ (where
$\beta_{\ell}=\rho_{\ell}^{-n/2}\|v\|_{L^{2}}(B_{\rho_{\ell}}(y_{\ell}))$), we can use Lemma~3.2 (the
local Schauder estimates) to first show that $v_{\ell}$ converges in $C^{1,\beta}$, for each
$\beta<\alpha$, locally in $B_{1}$ to give a $C^{1,\alpha}$ symmetric 2-valued harmonic function
$\varphi$ with $0\in{\cal{}K}_{\varphi}$, $\|\varphi\|_{L^{2}(B_{1})}\le 1$ and
$2^{n/2}\|\varphi\|_{L^{2}(B_{1/2})}\ge 2^{-3/2+\epsilon}$. However by Lemma~4.1 we know that
${\cal{}N}_{\varphi}(0)\ge 3/2$ and then by Remark~2.3(4) we have
$2^{n/2}\|\varphi\|_{L^{2}(B_{1/2})}\le 2^{-3/2}$ a contradiction. So as claimed we have
$(\rho/2)^{-n/2}\|v\|_{L^{2}(B_{\rho/2}(y))}\le
2^{-3/2+\epsilon}\rho^{-n/2}\|v\|_{L^{2}(B_{\rho}(y))}$ for each $\rho\le\rho_{0}$, and by iteration
(taking $\rho=2^{-j}\rho_{0}$, $j=1,2,\ldots,$) we conclude that we have the uniform decay
$\sup_{y\in B_{1/2}\cap {\cal{}K}_{u},\rho\le
\rho_{0}}\rho^{-3/2+\epsilon-n/2}\|v\|_{L^{2}(B_{\rho}(y))}<\infty$, and then the claimed estimate
$d(x)^{-\fr 1/2+\epsilon}|Dv(x)| <\infty$ by using elliptic estimates for $v$ (as $\pm$ a single valued
solution in balls which do not intersect ${\cal{}K}_{u}$).

Note that the above estimates for $v$ guarantee in particular that
$$%
\rho^{-1/2+\epsilon}[D_{i}v^{\lambda}D_{j}v^{\kappa}]_{\fr 1/2-\epsilon,B_{\rho}(y)} \le\beta, %
\quad y\in B_{1/2}\cap {\cal{}K}_{u}, \rho\in (0,\fr 1/8) %
\leqno{(1)}
$$%
for some constant $\beta$.

The argument giving the estimates for $u_{a}$ is similar, although since $u_{a}$ is single-valued we
can use the standard $C^{1,\alpha}$ theory rather than the $C^{1,\alpha}$ estimates of Lemma~3.2.
 In fact for each $y\in B_{1/2}$ and $\kappa=1,\ldots,k$ let
$$%
U_{y}^{\kappa}(x)=u^{\kappa}_{a}(x)-(x-y)\cdot Du_{a}^{\kappa}(y)-u_{a}^{\kappa}(y)
$$%
so that $U_{y}(y)=0$ and $DU_{y}(y)=0$ and 5.3', 5.12 and~(1) imply
$$%
D_{i}(A^{ij}D_{j}U_{y}^{\kappa}) = D_{j}F_{j} \text{ on }B_{1/2}, \, %
\rho^{-\fr 1/2+\epsilon}[F_{j}]_{\fr 1/2-\epsilon,B_{\rho}(y)}\le C\beta, %
\quad y\in B_{1/2}\cap {\cal{}K}_{u}, \rho\in (0,\fr 1/8), %
\leqno{(2)}
$$%
where $\beta$ is as in~(1), and so the standard $C^{1,\alpha}$ Schauder theory implies
\lal{\tag*{(3)} %
  \rho^{-\fr 1/2+\epsilon}[Du_{a}]_{\fr 1/2 -\epsilon,B_{\theta\rho}(y)} &\equiv %
  \rho^{-\fr 1/2+\epsilon}[DU_{y}]_{\fr 1/2 -\epsilon,B_{\theta\rho}(y)} \\ %
  &\hskip-.5in\le C\beta + C\rho^{-n/2-1}\|U_{y}\|_{L^{2}(B_{\rho}(y))}, %
  y\in B_{1/2}\cap {\cal{}K}_{u}, \rho\in (0,\fr 1/8), \quad C=C(\theta). %
} %
We claim there are $\rho_{0}=\rho_{0}(u)\in (0,\fr 1/8)$ and
$\lambda_{0}=\lambda_{0}(u)\ge 2$ such that
$$%
(\rho/2)^{-n/2-1}\|U_{y}\|_{L^{2}(B_{\rho/2}(y))} %
\le \max\{\rho^{-n/2-1}\|U_{y}\|_{L^{2}(B_{\rho}(y))},\lambda_{0}\beta\},\, %
y\in B_{1/2}\cap {\cal{}K}_{u}, \rho\in (0,\rho_{0}). %
\leqno{(4)}
$$%
If this fails then there are sequences $\rho_{\ell}\downarrow 0$, $\lambda_{\ell}\uparrow \infty$ and
$y_{\ell}\in B_{1/2}\cap {\cal{}K}_{u}$ with
$$%
\max\{\rho_{\ell}^{-n/2-1} %
\|U_{y_{\ell}}\|_{L^{2}(B_{\rho_{\ell}}(y_{\ell}))},\lambda_{\ell}\beta\} < %
(\rho_{\ell}/2)^{-n/2-1}\|U_{y_{\ell}}\|_{L^{2}(B_{\rho_{\ell}/2}(y_{\ell}))} %
\leqno{(5)}
$$%
and by~(3)
$$%
\rho_{\ell}^{-\fr 1/2+\epsilon}[DU_{y_{\ell}}]_{\fr 1/2 -\epsilon,B_{\theta\rho_{\ell}}(y_{\ell})} \le %
C(\theta) \max\{\rho_{\ell}^{-n/2-1} %
\|U_{y_{\ell}}\|_{L^{2}(B_{\rho_{\ell}}(y_{\ell}))},\beta\} %
\leqno{(6)}
$$%
for each $\theta\in (0,1)$. In particular~(5) implies
$$%
\beta\le \lambda_{\ell}^{-1} (\rho_{\ell}/2)^{-n/2-1} \|U_{y_{\ell}}\|_{L^{2}(B_{\rho_{\ell}/2}(y_{\ell}))} $$%
Defining $w_{\ell}(x)= \rho_{\ell}^{n/2}\|U_{y_{\ell}}\|_{L^{2}(B_{\rho_{\ell}}
(y_{\ell}))}^{-1}U_{y_{\ell}}(y_{\ell}+\rho_{\ell}x)$ for $x\in B_{1}$, we then have a subsequence of
$w_{\ell}$ which converges locally in $B_{1}$ with respect to the $C^{1,\gamma}$ norm for each
$\gamma<\ha-\epsilon$ to a harmonic function $\varphi$ with $\varphi(0)=0$, $D\varphi(0)=0$,
$\|\varphi\|_{L^{2}(B_{1})}\le 1$, $2^{n/2+1}\|\varphi\|_{L^{2}(B_{1/2})}\ge 1$ which is clearly
impossible (because, by the estimates of~2.3(4) with ${\cal{}N}_{\varphi}(0)\ge 2$,
$\|\varphi\|_{L^{2}(B_{1})}\ge \rho^{-n/2-2}\|\varphi\|_{L^{2}(B_{\rho})}$ for every $\rho<1$ in case
$\varphi$ is single-valued harmonic on $B_{1}$ with $\varphi(0)=0$ and $D\varphi(0)=0$). 
Thus~(4) is established, and can be written
\begin{align*}%
  &\max\{(\rho/2)^{-n/2-1}\|U_{y}\|_{L^{2}(B_{\rho/2}(y))},\lambda_{0}\beta\} \\ %
  &\hskip1in \le\max\{\rho^{-n/2-1}\|U_{y}\|_{L^{2}(B_{\rho}(y))},\lambda_{0}\beta\},\, %
  y\in B_{1/2}\cap {\cal{}K}_{u}, \rho\in (0,\rho_{0}). %
\end{align*}%
By an iteration similar to that used for $v$, we thus have
$$%
\sup_{y\in B_{1/2}\cap {\cal{}K}_{u},\rho\in (0,\rho_{0})} %
\rho^{-n/2-1}\|U_{y}\|_{L^{2}(B_{\rho}(y))}<\infty.
$$%
Then by~(3) we have the estimate
$$%
\sup_{y\in B_{1/2}\cap{\cal{}K}_{u},\rho\in (0,\rho_{0})} %
\rho^{-\fr 1/2+\epsilon}[Du_{a}]_{\fr 1/2 -\epsilon,B_{\theta\rho}(y)} <\infty %
$$%
which evidently implies
$$%
\sup_{y\in B_{1/2}\cap {\cal{}K}_{u},x\in B_{1/2}}|x-y|^{- 1+\epsilon} |Du_{a}(x)-Du_{a}(y)| <\infty,
$$%
and by using the non-divergence form~5.11' on balls $B_{\rho}(x)\subset
B_{1}\setminus{\cal{}K}_{u}$ where $x\in B_{1/2}\setminus{\cal{}K}_{u}$ and
$\rho<\min\{\rho_{0},d(x)/2\}$ we get the remaining estimate $d(x)^{\epsilon}|D^{2}u_{a}|\le C$ as
claimed.

\medskip

 \begin{state}{\bf{}5.14 Lemma.} %
   $D^{2}u_{a}\in L^{p}(B_{\rho}(y))$ for each $1\le p<\infty$ and $D^{2}v \in L^{2}(B_{\rho}(y))$
for each ball $B_{\rho}(y)$ with $\ovl{B}_{\rho}(y)\subset B_{1}$, and, for each $\theta\in (0,1)$,
$\|D^{2}u_{a}\|_{L^{2}(B_{\theta\rho}(y))}\le C\rho^{-2}\|u_{a}\|_{L^{2}(B_{\rho}(y))}$ and
$\|D^{2}v\|_{L^{2}(B_{\theta\rho}(y))}\le C\rho^{-2}\|v\|_{L^{2}(B_{\rho}(y))}$.
\end{state}%

{\bf{}Proof of Lemma~5.14:} 5.10 and the second identity in~5.3' implies
$$%
\sum_{} D_{i}\bigl(A^{ij}(Du_{a},Dv)D_{j}v^{\kappa}\bigr)= %
- \sum E^{ij\ell}_{\lambda}(Du_{a},Dv)D_{\ell}v^{\lambda} D_{i}D_{j}u_{a}^{\kappa}, %
$$%
and by using~5.13 to bound the right side, we see that this can be written
$$%
\sum_{} D_{i}\bigl(A^{ij}(Du_{a},Dv)D_{j}v^{\kappa}\bigr)= f^{\kappa},\quad |f^{\kappa}| \le
C\epsilon_{0}.
$$%
With $\zeta\in C_{c}^{\infty}(B_{\rho}(y))$, we multiply this equation by
$D_{\ell}(\gamma_{\delta}(D_{\ell}v^{\kappa})\zeta^{2})$, where
$\gamma_{\delta}(t)=\sgn(t)\max\{|t|-\delta,0\}$, and we integrate over $B_{\rho}(y)$. Integrating by
parts twice on the left then gives the identity
$$%
\int_{B_{\rho}(y)}D_{\ell}\bigl(A^{ij}(Du_{a},Dv)D_{j}v^{\kappa}\bigr) %
D_{i}\bigl(\gamma_{\delta}(D_{\ell}v^{\kappa})\zeta^{2}\bigr) = %
\int_{B_{\rho}(y)} f^{\kappa} D_{\ell}(\gamma_{\delta}(D_{\ell}v^{\kappa})\zeta^{2})
$$%
In view of~5.13 and 5.7 we see that
\begin{align*}%
  |D_{\ell}[A^{ij}(Du_{a},Dv)] | |Dv|&\le %
  ( |D_{p}A(Du_{a},Dv)||D^{2}u_{a}|+|D_{q}A(Du_{a},Dv)||D^{2}v|) |Dv| \\ %
  & \le C(|Du||D^{2}u_{a}| + |Dv||D^{2}v|)|Dv|\le C\epsilon_{0},
\end{align*}%
so this identity can be written
$$%
\int_{B_{\rho}(y)} A^{ij}(Du_{a},Dv)D_{\ell}D_{j}v^{\kappa} %
D_{i}\bigl(\gamma_{\delta}(D_{\ell}v^{\kappa})\zeta^{2}\bigr) = %
\int_{B_{\rho}(y)} (f^{\kappa}\delta_{i\ell} +g^{\kappa}) %
D_{i}(\gamma_{\delta}(D_{\ell}v^{\kappa})\zeta^{2})
$$%
with
$$%
|g^{\kappa}| \le C\epsilon_{0}.
$$%
Choosing $\zeta$ as a standard cut-off function in the ball $B_{\rho}(y)$ with $\zeta\equiv 1$ in
$B_{\rho/2}(y)$, and using the Cauchy-Schwarz inequality on the right, we conclude
$$%
\rho^{-n}\int_{B_{\rho/2}(y)\cap \{|D_{\ell}v^{\kappa}|>\delta\}} |DD_{\ell}v^{\kappa}|^{2}\le %
C\rho^{-2}\epsilon_{0}^{2}
$$%
with $C$ independent of $\delta$. Letting $\delta\downarrow 0$ and summing over $\ell,\kappa$
we thus obtain $D^{2}v\in L^{2}(B_{\rho/2}(y))$ with
$\rho^{-n/2}\|D^{2}v\|_{L^{2}(B_{\rho/2}(y))}\le C\rho^{-1}\epsilon_{0}$.

In view of Remark~5.12 we can now use a standard quasilinear elliptic difference quotient argument
on the equation 5.11 to establish that $D^{2}u_{a}\in L^{2}_{\text{loc}}(B_{1})$. 

To prove that $D^{2}u_{a}\in L^{p}_{\text{loc}}(B_{1})$ for $p\ge 2$ we first observe that, since
$D^{2}u_{a},D^{2}v\in L^{2}_{\text{loc}}(B_{1})$ by the above discussion, 5.3 now holds globally on
$B_{1}$ both pointwise a.e.\ and in the weak send, and so, by virtue of Remark~5.12, 5.11' holds
a.e.\ on $B_{1}$ and can be written 
$$%
\sum A^{ij}D_{i}D_{j}u^{\ell}_{a} = F^{\kappa},\quad |F^{\kappa}|\le C|Dv||D^{2}v|=C|D^{2}v|^{\fr
2/{p}}(|Dv||D^{2}v|^{1-\fr 2/{p}})\le C|D^{2}v|^{\fr 2/{p}}
$$%
by virtue of~5.13. Since $|D^{2}v|^{\fr 2/{p}}\in L^{p}_{\text{loc}}(B_{1})$ (because $|D^{2}u_{a}|\in
L^{2}_{\text{loc}}(B_{1})$) we can then use the standard interior Calderon-Zygmund $L^{p}$
estimates to prove $D^{2}u_{a}\in L^{p}_{\text{loc}}(B_{1})$ as claimed.
  
\section*{6\quad Some Growth Results for a Class of Linear Equations}

Here we consider $C^{1,\alpha}(B_{1},\R^{k})\cap W^{2,2}(B_{1},\R^{k})$ functions
$w=(w^{1},\ldots,w^{k})$ which are either single-valued or 2-valued symmetric, and we suppose
$w$ is a solution of a linear equation
$$%
D_{j}(a^{ij}_{\kappa\lambda}(x)D_{i}w^{\lambda}) =0 %
\text{ on }B_{1}\setminus {\cal{}K}_{w}, \quad \kappa=1,\ldots,k, %
\leqno{\hbox{\bf 6.1}}
$$%
where ${\cal{}K}_{w}=\{x:w(x)=0\text{ and }Dw(x)=0\}$ in the single-valued case and
${\cal{}K}_{w}=\{x:w(x)=\{0,0\}\text{ and }Dw(x)=\{0,0\}\}$ in the 2-valued case; of course in the
2-valued case we as usual take~6.1 to mean that we can locally, near each point of $B_{1}\setminus
{\cal{}K}_{w}$ write $w=\pm w_{1}$ with $w_{1}$ a $C^{1,\alpha}$ weak solution of the equation. 
We also assume
$$%
a^{ij}_{\kappa\lambda}= %
\delta_{ij}\delta_{\kappa\lambda}+b^{ij}_{\kappa\lambda},%
\quad b^{ij}_{\kappa\lambda}(0)=0,\quad %
[b^{ij}_{\kappa\lambda}]_{\alpha,B_{1}}\le \beta, %
\leqno{\hbox{\bf 6.2}}
$$%
and $w\in W^{2,2}(B_{1})$, with the estimates
$$%
\|D^{2}w\|_{L^{2}(B_{\rho/2})}\le \beta\rho^{-2}\|w\|_{L^{2}(B_{\rho})}, \quad \rho\in (0,1),
\leqno{\hbox{\bf 6.3}}
$$%
where $\beta>0$, $\alpha\in (0,1)$ in the single-valued case, and $\alpha\in (0,\delta)$ with
$\delta=\delta(n)$ as in~3.2 in the 2-valued case.  

\medskip

Observe that, subject only to 6.1, 6.2 and the assumption $w\in C^{1,\alpha}(B_{1},\R^{k})$, the
Schauder estimates in Lemma~3.2 are applicable to give
$$%
\rho^{\alpha}[Dw]_{\alpha,B_{3\rho/4}}+ %
                               \sup_{B_{3\rho/4}}|Dw| + \rho^{-1}  \sup_{B_{3\rho/4}}|w|\le %
C(\rho^{-2-n}\int_{B_{\rho}}|w|^{2})^{1/2} %
\leqno{\hbox{\bf 6.4}}
$$%
for $\rho<\rho_{0}$, where $\rho_{0}=\rho_{0}(n,k,\beta)\in (0,1)$ is suitably small.

Now, for the moment, assume $|w(0)|=0$, let $\rho\in (0,1),\theta\in (0,1/2)$ and that $\lambda$ is
the mean value of $|w|^{2}$ on $B_{\theta\rho}$ (so that
$\lambda=|B_{\theta\rho}|^{-1}\int_{B_{\theta\rho}}|w|^{2}$). Then the appropriate version of the
Poincar\` e inequality (\cite[\S7.8]{GilT83}) on $B_{\rho}$ implies
$$%
\sint_{B_{\rho}}||w|^{2}-\lambda| \le \gamma\rho\sint_{B_{\rho}}|D|w|^{2}|, %
\quad \gamma=\gamma(n,\theta),
$$%
whence
$$%
\sint_{B_{\rho}}|w|^{2}\le\lambda|B_{\rho}| + \gamma\rho\sint_{B_{\rho}}|D|w|^{2}|,
$$%
and on the other hand $\lambda=|w|^{2}(y)$ for some $y\in B_{\theta\rho}$ and $|w|^{2}(0)=0$, so
by 1-variable calculus we have $\lambda\le \theta\rho\sup_{B_{\theta\rho}}|D|w|^{2}|\le
2\theta\rho\sup_{B_{\theta\rho}}(|Dw||w|)$ and hence by~6.4
$$%
\sint_{B_{\rho}}|w|^{2}\le C\theta \sint_{B_{\rho}}|w|^{2} +
\gamma\rho\sint_{B_{\rho}}|D|w|^{2}|,\quad C=C(n,k,\beta),
$$%
hence by choosing $\theta=(2C)^{-1}$ and using Cauchy-Schwarz we get
$$%
\sint_{B_{\rho}}|w|^{2}\le C\rho^{2}\sint_{B_{\rho}}|Dw|^{2}, \quad C=C(n,k,\beta),
\leqno{\hbox{\bf 6.5}}
$$%
provided 6.1, 6.2 hold and $|w(0)|=0$.

\bigskip

We next show that a ``doubling condition'' for $w$ is sufficient to establish a suitable monotonicity
and bounds for the frequency function $N_{w}(\rho)$ of $w$, which is the function of $\rho\in
(0,1)$ defined by
$$%
N_{w}(\rho)= %
\Fr{\rho^{2-n}\mint_{B_{\rho}}|Dw|^{2}}/{\rho^{1-n}\mint_{\partial B_{\rho}}|w|^{2}}
$$%
{\abovedisplayskip1pt\belowdisplayskip1pt
 \begin{state}{\bf{}6.6 Lemma.} %
   If $w\in C^{1,\alpha}(B_{1})$, if the hypotheses~6.1, 6.2, 6.3 hold, if $0\in{\cal{}K}_{w}$, and if
there is $\gamma>1$ and $\sigma\in (0,1/2]$ such that
$$%
\|w\|_{\rho}\le \gamma \|w\|_{\rho/2} \quad \forall \,\rho\in (0,\sigma],
\hskip1em\text{where}\hskip5pt \|w\|_{\rho}=\Bigl(\rho^{1-n}\sint_{\partial
B_{\rho}}w^{2}\Bigr)^{1/2},
$$%
\vskip-10pt then
$$%
\sint_{B_{\rho}}w^{2}\le C\gamma^{2} \sint_{B_{\rho/2}}w^{2} \quad \forall \,\rho\in (0,\sigma], %
\quad C=C(n),
$$%
\vskip-10pt and
$$%
{\Fr d/{d\rho}}\Bl(e^{\Lambda \rho^{\alpha}}N_{w}(\rho)\Br)\ge 0 \quad \forall\,\rho\in (0,\sigma],
\quad\text{ where } \Lambda=\Lambda(n,\alpha,\beta,\gamma).
$$%
Furthermore, for each sequence $\rho_{j}\downarrow 0$ there is a subsequence $\rho_{j^{\prime}}$ such that \newline %
$\rho_{j'}^{n/2}\|w\|_{L^{2}(B_{\rho_{j^{\prime}}})}^{-1}w(\rho_{j^{\prime}}x)\to \varphi(x)$ in
$C^{1}$ locally on $\R^{n }$, where $\varphi$ is a homogeneous harmonic 2-valued symmetric
$\smash{C^{1,\alpha}}$ in the two-valued case, and smooth harmonic in the single-valued case, and
\begin{align*}%
  \noalign{\vskip-1pt} &e^{\Lambda \rho^{\alpha}}N_{w}(\rho)\ge %
  {\cal{}N}_{w}(0)={\cal{}N}_{\varphi}(0)\ge 3/2 %
  \quad \forall\,\rho\in (0,\sigma] \text{ in the 2-valued case,} \\ %
  \noalign{\vskip-3pt} %
  &\hskip3in\text{hence }\|w\|_{\rho} \le %
  C\|w\|_{\sigma}(\rho/\sigma)^{3/2} \text{ for }\rho\in (0,\sigma] \\ %
  &e^{\Lambda \rho^{\alpha}}N_{w}(\rho)\ge %
  {\cal{}N}_{w}(0)={\cal{}N}_{\varphi}(0)\ge 2\quad \forall\,\rho\in (0,\sigma] \text{ in the single-valued case,} \\ %
  \noalign{\vskip-3pt} %
  &\hskip3in\text{hence }\|w\|_{\rho}\le C\|w\|_{\sigma}(\rho/\sigma)^{2}\text{ for }\rho\in
(0,\sigma],
\end{align*}%
\vskip-5pt where $C=C(n,\alpha,\beta,\gamma)$.
\end{state}}%

{\textbf{Remark:}} Note that $\Lambda,C$ do not depend on $\sigma$ in the above lemma.

\smallskip
            
{\textbf{Proof of Lemma 6.6:}} First observe that by integrating the inequality
$$%
\sint_{\partial B_{\rho}}w^{2}\le 2^{n-1}\gamma^{2} \sint_{\partial B_{\rho/2}}w^{2} %
\quad \forall \,\rho\in (0,\sigma] %
\leqno{\hbox{(1)}}
$$%
with respect to $\rho$ we conclude immediately that
$$%
\sint_{B_{\rho}}w^{2}\le C\gamma^{2} \sint_{B_{\rho/2}}w^{2} %
\quad \text{ for all $\rho\in (0,\sigma]$} %
\leqno{\hbox{(2)}}
$$%
as claimed.  Define $N_{w}$ as above, i.e.\
$$%
N_{w}(\rho)= %
\Fr{\rho^{2-n}\mint_{B_{\rho}}|Dw|^{2}}/{\rho^{1-n} %
  \mint_{\partial B_{\rho}}|w|^{2}}. %
\leqno{\hbox{(3)}}
$$%
Since $w^{\lambda}D_{j}w^{\kappa}$ and $D_{i}w^{\lambda}D_{j}w^{\kappa}$ are single valued
$W^{1,1}$ functions on $B_{1}$ with $D_{\ell}(w^{\lambda}D_{j}w^{\kappa})
=(D_{\ell}w^{\lambda})D_{j}w^{\kappa}+ w^{\lambda}D_{\ell}D_{j}w^{\kappa}$ and
$D_{\ell}(D_{i}w^{\lambda}D_{j}w^{\kappa})=(D_{\ell}D_{i}w^{\lambda})
D_{j}w^{\kappa}+D_{i}w^{\lambda}D_{\ell}D_{j}w^{\kappa}$ a.e.\ on $B_{1}$, and since we can
use the weak form of~6.1 to check that $\int_{B_{1}}(\Delta w^{\kappa})w^{\kappa}\zeta=
-\int_{B_{1}}b^{i\ell}_{\kappa\lambda}w^{\lambda}_{i}D_{\ell}(w^{\kappa}\zeta)$ and
$\int_{B_{1}}(\Delta w^{\kappa})w^{\kappa}_{j}\zeta_{j}=
-\int_{B_{1}}b^{i\ell}_{\kappa\lambda}w^{\lambda}_{i}D_{\ell}(w^{\kappa}_{j}\zeta_{j})$ for
$\zeta,\zeta_{j}\in C^{1}_{c}(B_{1})$, it is straightforward to check the two identities 
$$%
\sint_{B_{1}}w^{\kappa}_{i}D_{i}\bigl(w^{\kappa}\zeta\bigr)=   %
\sint_{B_{1}} b^{ij}_{\kappa\lambda}w^{\lambda}_{i}D_{j}\bigl(w^{\kappa}\zeta\bigr), %
\leqno{\hbox{(4)}}
$$%
and
$$%
\sint_{B_{1}} \bigl(|Dw^{\kappa}|^{2} %
\delta_{ij} - 2w_{i}^{\kappa}w_{j}^{\kappa}\bigr)%
D_{i}\zeta_{j} = %
2\sint_{B_{1}}b^{i\ell}_{\kappa\lambda}w_{i}^{\lambda}(w_{j}^{\kappa}D_{\ell}\zeta_{j}+  %
                                                                                       D_{\ell}w_{j}^{\kappa}\zeta_{j}) %
\leqno{\hbox{(5)}}
$$%
for any Lipschitz functions $\zeta,\zeta_{1},\ldots,\zeta_{n}$ with compact support in $B_{1}$,
where $w^{\lambda}_{i}=D_{i}w^{\lambda}$ and repeated indices indicate summation as usual. 
Notice that the second identity is checked directly by integrating by parts in the expression on the
left, using the fact that $D^{2}w\in L^{2}$ and $D_{j}(a^{\ell
j}_{\kappa\lambda}D_{\ell}w^{\lambda})=0$ a.e.\ on $B_{1}$.  Using~6.2 to give
$|b^{ij}_{\kappa\lambda}|\le \beta\rho^{\alpha}$ in each of these identities, and letting $\zeta$
approximate the indicator function of $B_{\rho}$ in the first identity and letting $\zeta_{i}$
approximate $x_{i}$ times in the indicator function of $B_{\rho}$ in the second identity, we obtain
the two key inequalities
$$%
\Bigl|\sint_{B_{\rho}}|Dw|^{2}-\sint_{\partial B_{\rho}}ww_{r}\Bigr|\le %
C\rho^{\alpha}\sint_{B_{\rho}}|Dw|^{2}+ C\rho^{\alpha}\sint_{\partial B_{\rho}}|w||Dw| %
\leqno{\hbox{(6)}}
$$%
and
\begin{align*}%
  \tag*{(7)} %
  &\Bigl|(n-2)\sint_{B_{\rho}}|Dw|^{2}-\rho\sint_{\partial B_{\rho}}|Dw|^{2}+ %
  2\rho\sint_{\partial B_{\rho}}w_{r}^{2}\Bigr| \\ %
  &\hskip.4in \le C\rho^{\alpha}\sint_{B_{\rho}}|Dw|^{2} + %
  C\rho^{\alpha+1}\sint_{\partial B_{\rho}}|Dw|^{2}+ %
  C\rho^{\alpha+1}\sint_{B_{\rho}}|Dw||D^{2}w|
 \end{align*}%
 Next we observe that by direct computation we have, weakly in $B_{1}$,
$$%
\Delta |w|^{2}= 2 |Dw|^{2} + 2b^{ij}_{\kappa\lambda}w^{\kappa}_{i}w^{\lambda}_{j}- 2D_{j}(
b^{ij}_{\kappa\lambda}w^{\lambda}w^{\kappa}_{i}).
$$%
Integrating this over $B_{\rho}$ and using the identity $\int_{B_{\rho}}\Delta
f=\rho^{n-1}{d\over{}d\rho}(\rho^{1-n}\int_{\partial B_{\rho}}f)$ with $f=|w|^{2}$, we then conclude
$$%
\rho^{n-1}\tfrac{d}{d\rho}(\rho^{1-n}\sint_{\partial B_{\rho}}|w|^{2}) = %
2\int_{B_{\rho}}(|Dw|^{2}+b_{\kappa\lambda}^{ij}w^{\kappa}_{i}w^{\lambda}_{j}) %
-2\rho^{-1}\int_{\partial B_{\rho}}x^{j}b_{\kappa\lambda}^{ij}w^{\kappa}_{i}w^{\lambda}, %
$$%
hence
$$%
\tfrac{d}{d\rho}(\sint_{\partial B_{\rho}}|w|^{2}) =(n-1)\rho^{-1}\sint_{\partial B_{\rho}}|w|^{2} %
+2\int_{B_{\rho}}(|Dw|^{2}+b_{\kappa\lambda}^{ij}w^{\kappa}_{i}w^{\lambda}_{j}) %
-2\rho^{-1}\int_{\partial B_{\rho}}x^{j}b_{\kappa\lambda}^{ij}w^{\kappa}_{i}w^{\lambda}, %
\leqno{(8)}
$$%
which by~6.2 evidently implies, for $\rho\le \rho_{0}(n,k,\beta)$,
$$%
\tfrac{d}{d\rho}\sint_{\partial B_{\rho}}|w|^{2} \ge -C\rho^{\alpha} \sint_{\partial B_{\rho}}|w||Dw|.
$$%
By integrating this with respect to $\rho$ we obtain
$$%
\sint_{\partial B_{\tau}}|w|^{2} \le \sint_{\partial B_{\rho}}|w|^{2} +
C\rho^{\alpha}\int_{B_{\rho}}|w||Dw|, \quad \tau<\rho\le \rho_{0},
$$%
where $\rho_{0}\in (0,1)$ depends only on $n,\beta$ and $k$. 
Integration with respect to $\tau\in (0,\rho)$ then shows
$$%
\sint_{B_{\rho}}|w|^{2} \le %
\rho\sint_{\partial B_{\rho}}|w|^{2} + C\rho^{\alpha+1}\int_{B_{\rho}}|w||Dw|, %
\quad \rho\le \rho_{0}. %
\leqno{(9)}
$$%
We can also use~(8) and the estimates~6.4 to give the upper bound
$$%
\tfrac{d}{d\tau}(\sint_{\partial B_{\tau}}|w|^{2})\le C\rho^{-2}\int_{B_{\rho}}|w|^{2}, \quad \tau \le
\rho/2,\, \rho <1,
$$%
which integrates to give
$$%
\sint_{\partial B_{\rho}}|w|^{2}\le C\rho^{-1}\int_{B_{2\rho}}|w|^{2},\quad \rho<1/2.  \leqno{(10)}
$$%
Thus if we have
$$%
\|w\|_{\rho}\le\gamma\|w\|_{\rho/2} %
\quad \rho\in(0,\sigma] \leqno{(11)}
$$%
then by~6.4, (8) and~(9) we have
$$%
\sint_{B_{\rho}}(\rho^{2}|Dw|^{2}+|w|^{2})+ %
\rho\sint_{\partial B_{\rho}} \rho^{2}|Dw|^{2} \le %
C\min\{\rho\sint_{\partial B_{\rho}}|w|^{2} ,\sint_{B_{\rho}}|w|^{2}\}, %
\quad \rho\in (0,\sigma], %
\leqno{\hbox{(12)}}
$$%
and, by~(10), (11), 6.4 and 6.5
$$%
\sint_{\partial B_{\rho}} |w|^{2} \le %
C\min\{\rho\sint_{B_{\rho}}|Dw|^{2} ,\rho^{-1}\sint_{B_{\rho}}|w|^{2}\} %
\leqno{\hbox{(13)}}
$$%
for each $\rho\in (0,\sigma]$.  Now let
$$%
{\cal{}D}(\rho)=\rho^{2-n}\int_{B_{\rho}}|Dw|^{2},\quad %
{\cal{}H}(\rho)=\rho^{1-n}\int_{\partial B_{\rho}}w^{2},
$$%
and observe that ${\cal{H}}(\rho)$ never vanishes for $\rho\in (0,\sigma]$ by virtue of~(11)
and~(12), so we can define
$$%
N(\rho)={\cal{}D}(\rho)/{\cal{}H}(\rho),
$$%
and note that by~(12) and~(13) we have
$$%
C^{-1} \le N(\rho)\le C, \quad \rho\in (0,\sigma].  \leqno{(14)}
$$%
Also by~(7), (12) and 6.3 we have
$$%
|(n-2)\sint_{B_{\rho}}|Dw|^{2}-\sint_{\partial B_{\rho}}|Dw|^{2}+ %
2\sint_{\partial B_{\rho}}w_{r}^{2}|\le %
C\rho^{\alpha}\min\{\rho\sint_{\partial B_{\rho}}|w|^{2} , %
\sint_{B_{\rho}}|w|^{2},\rho^{2}\sint_{B_{\rho}}|Dw|^{2} \} %
\leqno{\hbox{(15)}}
$$%
for every $\rho\in (0,\sigma]$, assuming~(11), which implies
$$%
{\cal{}D}^{\prime}(\rho)\ge 2\sint_{\partial B_{\rho}}|w_{r}|^{2} - %
C\rho^{\alpha}\min\{\rho\sint_{\partial B_{\rho}}|w|^{2}, %
\sint_{B_{\rho}}|w|^{2},\rho^{2}\sint_{B_{\rho}}|Dw|^{2} \} %
\quad \forall \,\rho\in (0,\sigma], %
\leqno{\hbox{(16)}}
$$%
Also, by virtue of~(6) we have
$$%
{\cal{}H}^{\prime}(\rho)=2\rho^{1-n}\sint_{\partial B_{\rho}}ww_{r}= %
\rho^{-1}{\cal{}D}(\rho)+E,\quad %
\text{ with } |E|\le \rho^{\alpha+1-n}\sint_{B_{\rho}}(\rho^{-2}|w|^{2}+|Dw|^{2}) %
\leqno{\hbox{(17)}}
$$%
Now
$$%
N^{\prime}(\rho)= %
\Fr{{\cal{}D}^{\prime}(\rho){\cal{}H}(\rho)-{\cal{}H}^{\prime}(\rho){\cal{}D}(\rho)}/ %
{{\cal{}D}(\rho){\cal{}H}(\rho)} %
$$%
and using~(15), (16) and (17) we then get
$$%
N^{\prime}(\rho) \ge \Fr{2{\cal{}H}(\rho)\mint_{\partial B_{\rho}}w_{r}^{2} - %
  2({\cal{}H}^{\prime}(\rho))^{2} - %
  C\rho^{\alpha-1} {\cal{}H}^{2}(\rho)}/{{\cal{}D}(\rho){\cal{}H}(\rho)}\ge %
-C\Fr{\rho^{\alpha-1}{\cal{}H}^{2}(\rho)}/{{\cal{}D}(\rho){\cal{}H}(\rho)}, %
$$%
where we used ${\cal{}H}(\rho)\int_{\partial B_{\rho}}w_{r}^{2} -
({\cal{}H}^{\prime})^{2}=\rho^{2-2n}\bigl(\int_{\partial B_{\rho}}w^{2}\int_{\partial
B_{\rho}}w_{r}^{2} - (\int_{\partial B_{\rho}}ww_{r})^{2}\bigr) \ge 0$ by Cauchy-Schwarz. Finally in
view of~(14) we thus have
$$%
N^{\prime}(\rho) \ge -C\rho^{\alpha-1}N(\rho),
$$%
which can be written
$$%
\Fr{d}/{d\rho}\Bigl(\exp(\alpha^{-1}C\rho^{\alpha})\,N(\rho)\Bigr)\ge 0 \leqno{(18)}
$$%
as claimed.

Now take a sequence $\rho_{j}\downarrow 0$ and let $w_{j}(x)=\lambda_{j}w(\rho_{j}x)$ with
$\lambda_{j}$ chosen to ensure that $\|w_{j}\|_{L^{2}(B_{1})}=1$ for each $j$.  By virtue of the
estimates~6.4 and~(2) we have $|w_{j}|_{C^{1,\alpha}(B_{R})}\le C(R,n,\gamma)$ for all $R>0$ and
sufficiently large $j$ depending on $R$, so a subsequence $w_{j^{\prime}}$ converges locally in
$\R^{n}$ in $C^{1}$ to a symmetric 2-valued $C^{1,\alpha}$ function $\varphi$ and $\varphi$ is
evidently harmonic. Furthermore $N_{w_{j^{\prime}}}(\rho)\to N_{\varphi}(\rho)$ for each fixed
$\rho$, whereas by the monotonicity~(18) we have $N_{w_{j^{\prime}}}(\rho)\to
\lim_{\sigma\downarrow 0} N_{w}(\sigma)={\cal{}N}_{w}(0)$, independent of $\rho$. So
$N_{\varphi}(\rho)\equiv {\cal{}N}_{w}(0)$ and hence $\varphi$ is homogeneous by 2.4(1), and
by~4.1 the order of homogeneity is $\ge 3/2$, so ${\cal{}N}_{\varphi}(0)\ge 3/2$ as claimed.

The rest of the proof, dealing with the single valued case, is similar.
  
\medskip

Next we prove a growth lemma for $w$:

\begin{state}{\bf{}6.7 Lemma.} %
  If $\alpha\in (0,\delta)$, $\delta=\delta(n)$ as in~3.2, there is
$\rho_{0}=\rho_{0}(n,k,\alpha,\beta)\in (0,1/2]$ and $C=C(n,\alpha,\beta)$ and
$\gamma=\gamma(n,k,\alpha,\beta)$ where $\gamma\in (3/2,3/2+\delta)$ in the 2-valued case and
$\gamma\in (2,3)$ in the single valued case, such that if the hypotheses~6.1, 6.2 hold and if we
write $\|w\|_{\rho}= (\rho^{1-n}\int_{\partial B_{\rho}}w^{2})^{1/2}$ for each $\rho\in (0,\rho_{0}]$,
then {\abovedisplayskip5pt\belowdisplayskip5pt
$$%
\|w\|_{\rho/2}\ge 2^{\gamma}\|w\|_{\rho/4} \Rightarrow \|w\|_{\sigma}\ge %
(2\sigma/\rho)^{\gamma} \|w\|_{\rho/2} \text{ for each }\sigma\in [3\rho/4,\rho]  %
\text{ and each }\rho\in(0, \rho_{0}], %
$$}%
so in particular (taking $\sigma=\rho$)
$$%
{\abovedisplayskip1pt\belowdisplayskip1pt %
  \|w\|_{\rho/2}\ge 2^{\gamma}\|w\|_{\rho/4} %
  \Rightarrow \|w\|_{\rho}\ge 2^{\gamma} \|w\|_{\rho/2} } %
$$%
for each $\rho\in (0,1]$.
\end{state}%

{\bf 6.8 Remark:} In particular this lemma shows that if $\sigma\in (0,\fr {\rho_{0}}/2]$ is such that
$\|w\|_{2\sigma}\ge 2^{\gamma}\|w\|_{\sigma}$ then $\|w\|_{2^{j}\sigma}\ge 2^{\gamma}
\|w\|_{2^{j-1}\sigma}$ for every $j=1,\ldots$ such that $2^{j}\sigma\le \rho_{0}$, which implies
$$%
\|w\|_{\rho}\le C\rho^{\gamma}\|w\|_{\rho_{0}}
$$%
for all $\rho\in (\sigma,\fr {\rho_{0}}/2]$, whence
$$%
\|w\|_{\rho}\le C\rho^{\gamma}\|w\|_{L^{2}(B_{1})} \text{ for all }\rho\in (\sigma,\fr 1/2].
$$%

\medskip
                      
{\textbf{Proof of Lemma~6.7:}} We give the proof first in the 2-valued case. If the contrary holds
with $\rho_{0}=1/\ell$, $\rho=\rho_{\ell}$ and $w=w_{\ell}$ for $\ell=1,2,\ldots$, then by rescaling
$x\to \rho^{-1}x$ and using the $C^{1,\alpha}$ estimates of~6.4 we see that a subsequence has a
2-valued symmetric $C^{1,\alpha}$ limit $\varphi$ on $B_{\sigma}$ for some $\sigma\in [3/4,1]$
which is harmonic on $B_{\sigma}\setminus{\cal{}K}_{\varphi}$ and with
$\rho^{-n/2-\gamma}\|\varphi\|_{L^{2}(\partial B_{\rho})}$, as a function of $\rho\in [1/4,\sigma]$,
taking a local maximum value at some $\rho_{0}\in (1/4,\sigma)$. By the Remark~2.3(1) we know
that $\log\|\varphi\|_{\rho}$ a convex function of $t=\log \rho$ for $\rho\in (1/4,\sigma)$ and hence
so is $\log(\rho^{-\lambda}\|\varphi\|_{\rho})= -\lambda t + \log\|\varphi\|_{\rho}$. A convex
function attaining a local interior maximum is constant, hence $\rho^{-\lambda}\|\varphi\|_{\rho}$
is constant for $\rho\in (1/4,\sigma)$. Using the Remark~2.4 we then have that $\varphi_{r}$ is a
constant multiple of $\varphi$ so that $\varphi$ extends to all of $\R^{n}$ as a homogeneous degree
$\gamma$ function. Since we take $\gamma\in (3/2,3/2+\delta)$ with $\delta$ as in Lemma~4.2,
this contradicts the result of Lemma~4.2.

The proof in case $w$ is single valued is similar, except that we use $\gamma\in (2,3)$ and since
there are no single valued harmonic functions on $\R^{n}$ which are homogeneous of non-integer
degree, this again gives a contradiction.

\medskip

Using Lemma~6.6 and Lemma~6.7, we can now prove the following regularity/decay result:

\begin{state}{\bf{}6.9 Theorem.} %
  If $\alpha\in (0,\delta)$, $\delta=\delta(n)$ as in~3.2, there is
$\epsilon=\epsilon(n,p,\alpha,\beta)\in (0,1/2]$ and $C=C(n,p,\alpha,\beta)$ such that if\/ $w\in
C^{1,\alpha}$, if 6.1, 6.2, 6.3 hold and if\/ $0\in{\cal{}K}_{w}$, then in the 2-valued symmetric case
we have {\abovedisplayskip3pt\belowdisplayskip3pt
$$%
\sup_{B_{\rho}}|w|\le C\|w\|_{L^{2}(B_{1})}\rho^{3/2},\quad \rho\in (0,\ha]
$$}%
and in the single-valued case we have {\abovedisplayskip3pt\belowdisplayskip3pt
$$%
\sup_{B_{\rho}}|w|\le C\|w\|_{L^{2}(B_{1})}\rho^{2},\quad \rho\in (0,\ha].
$$}%
\end{state}%
{\bf{}Proof:} We give the proof in the 2-valued case first: First assume $0\in {\cal{}K}_{w}\cap
B_{1/2}$. Let $\delta=\delta(n)\in (0,1)$ be as in the 2-valued case of Lemma~6.6, and let
$\gamma=3/2+\delta/2$. For suitable $\epsilon=\epsilon(n,p)>0$, we can apply the 2-valued case of
Lemma~6.6. So let $\sigma\in [0,1]$ be $\inf\{\ha,\{\rho\in (0,1/2]: \|w\|_{\rho}\ge 2^{\gamma}
\|w\|_{\rho/2}\}\}$. Then in accordance with Remark~6.8 we have
$$%
\|w\|_{\rho}\le C\rho^{\gamma}\|w\|_{L^{2}(B_{1})}\le %
C\rho^{3/2}\|w\|_{L^{2}(B_{1})}, \quad \rho\in (\sigma,1/2] %
\leqno{\hbox{(1)}}
$$%
and, assuming $\sigma\neq 0$, $\|w\|_{\rho}< 2^{\gamma}\|w\|_{\rho/2}$ for every $\rho\in
(0,\sigma]$, and hence by the 2-valued case of Lemma~6.6 we have
$$%
\|w\|_{\rho}\le C\|w\|_{\sigma}(\rho/\sigma)^{3/2},\quad \rho\in (0,\sigma]. %
\leqno{\hbox{(2)}}
$$%
Thus by combining~(1) and~(2) we have
$$%
\|w\|_{\rho}\le C\|w\|_{L^{2}(B_{1})}\rho^{3/2}, \quad \forall\, \rho\in (0,1/2], %
$$%
and by integration with respect to $\rho$
$$%
\Bl(\rho^{-n}\sint_{B_{\rho}}w^{2}\Br)^{1/2} \le C\|w\|_{L^{2}(B_{1})}\rho^{3/2}, %
\quad\forall\, \rho\in (0,1/2].
$$%
Next observe that by the $C^{1,\alpha}$ estimates we then have
$$%
\sup_{B_{\rho/2}}|w| \le C\Bl(\rho^{-n}\sint_{B_{\rho}}w^{2}\Br)^{1/2} \le %
C\|w\|_{L^{2}(B_{1})}\rho^{3/2},\quad\forall\, \rho\in (0,1/2]. %
$$%
This completes the proof in the 2-valued symmetric case.  The proof in the single-valued case is
similar except that we use $\gamma\in (2,3)$ and the single-valued cases of Lemmas~6.6, 6.7, and
in place of~4.4 we use the standard fact that each single-valued homogeneous harmonic function on
$\R^{n}$ is given by a homogeneous harmonic polynomial and hence has integer order of
homogeneity.

 \section*{7\quad Regularity for  %
           \bm{$u_{a}=\fr 1/2(u_{1}+u_{2})$}    and  %
                          \bm{$v=\pm\ha(u_{1}-u_{2})$}, Part~II}  %
                        We first observe that, in view of Lemmas~5.13 and~5.14, the equation~5.10 has the
correct form to ensure that Theorem~6.9 can be applied with $w=v$, whence we have

\begin{state}{\bf{}7.1 Theorem.} %
  $v=\{\pm\fr 1/2(u_{1}-u_{2})\}$ is locally $C^{1,1/2}$ in $B_{1}$, and we have the estimates
$$%
|v(x)| \le C\epsilon_{0}d(x)^{3/2},\,\,\, |Dv(x)|\le C\epsilon_{0} d(x)^{1/2}, %
\,\,\, |D^{2}v(x)|\le C \epsilon_{0} d(x)^{-1/2}, \quad x\in B_{1/2}, %
$$%
where $C=C(n,k)$ and $d(x)=\dist(x,{\cal{}K}_{u})$.
\end{state}%

In view of~5.11 we can also apply the single-valued case of Theorem~6.9 to the average
$u_{a}=\ha(u_{1}+u_{2})$ in balls centered at $0$ (since we assume $u(0)=\{0,0\}$ and
$Du(x)=\{0,0\}$). Thus we have
$$%
\sup_{B_{\rho}}|u_{a}|\le C\epsilon_{0}\rho^{2} \quad \forall \,\rho\in (0,1].  \leqno{\text{\bf 7.2}}
$$%
We want to show that similar decay estimates hold for $u_{a}(x)-u_{a}(x_{0})-(x-x_{0})\cdot
Du_{a}(x_{0})$ on balls $B_{\rho}(x_{0})$ for any $x_{0}\in {\cal{}K}_{u}\cap B_{1/4}$:

\begin{state}{\bf{}7.3 Lemma.} %
  If $x_{0}\in B_{1/4}\cap {\cal{}K}_{u}$ then
$$%
|u_{a}(x)-u_{a}(x_{0})-(x-x_{0})Du_{a}(x_{0})| \le C\epsilon_{0}\rho^{2} %
$$%
for all $x\in B_{\rho}(x_{0}), \, \,\rho\in (0,1/4]$.
\end{state}%
{\bf{}Proof:}
Since $|Du_{a}|<\epsilon_{0}$, we can choose an $(n+k)\times (n+k)$ orthogonal matrix ${\cal{}Q}$
with $|I-{\cal{}Q}|<C\epsilon_{0}$ such that ${\cal{}Q}$ takes the tangent space of $\graph u_{a}$ at
at the point $(x_{0},u(x_{0}))$ to the space $\R^{n}\times\{0\}$ and $\graph u=\graph
\{u_{1},u_{2}\}$ is transformed to $\graph \{\wtilde u_{1},\wtilde u_{2}\}$, where $\wtilde
u=\{\wtilde u_{1},\wtilde u_{2}\}$ is a $C^{1,1/2}$ function over $B_{1/4}$ with $[D\wtilde
u]_{1/2}<C\epsilon_{0}$,
\begin{align*}%
  \tag*{(1)}(\xi,\wtilde u_{1}(\xi))&=(x-x_{0}, u_{1}(x)-u_{a}(x_{0})){\cal{}Q} \\
  (\eta,\wtilde u_{2}(\eta))&=(x-x_{0},u_{2}(x)-u_{a}(x_{0})){\cal{}Q}.
\end{align*}%
Notice that since $(e_{i},D_{i}u_{a}(x_{0})$ is in the tangent space of $\graph u_{a}$ at
$(x_{0},u_{a}(x_{0}))$ for $i=1,\ldots,n$ we have $(e_{i},D_{i}u_{a}(x_{0}))Q_{k}=0$, where $Q_{k}$
is the $(n+k)\times k$ matrix consisting of the last $k$ columns of ${\cal{}Q}$, so for $i=1,\ldots,n$
we have
$$%
(e_{i},0)Q_{k}=-D_{i}u_{a}(x_{0})Q_{kk}, %
\leqno{(2)}
$$%
where $Q_{kk}$ is the $k\times k$ matrix consisting of the last $k$ rows of $Q_{k}$. On the other
hand since $x-x_{0}=\sum_{i=1}^{n}(x_{i}-x_{0i})e_{i}$ we see from~(1) that
\begin{align*}%
  &\tsum_{i=1}^{n}(x_{i}-x_{0i})(e_{i},0)Q_{k}+(u_{1}(x)-u_{a}(x_{0}))Q_{kk}=\wtilde u_{1}(\xi) \\
  &\tsum_{i=1}^{n}(x_{i}-x_{0i})(e_{i},0)Q_{k}+(u_{2}(x)-u_{a}(x_{0}))Q_{kk}=\wtilde u_{2}(\eta),
\end{align*}%
and so by~(2)
\begin{align*}%
  &\bigl(-(x-x_{0})Du_{a}(x_{0})+(u_{1}(x)-u_{a}(x_{0}))\bigr)Q_{kk}=\wtilde u_{1}(\xi) \\
  &\bigl(-(x-x_{0})Du_{a}(x_{0})+(u_{2}(x)-u_{a}(x_{0}))\bigr)Q_{kk}=\wtilde u_{2}(\eta),
\end{align*}%
and by taking sums we have
\begin{align*}%
  &\bigl(-(x-x_{0})Du_{a}(x_{0})+(u_{a}(x)-u_{a}(x_{0}))\bigr)Q_{kk}= %
  \ha (\wtilde u_{2}(\eta)+\wtilde u_{1}(\xi)) \\ %
  &\hskip.8in %
  =\ha(\wtilde u_{a}(\xi)+\wtilde u_{a}(\eta))+\tfrac{1}{4} ((\wtilde u_{1}(\xi)-\wtilde u_{2}(\xi)) %
  -(\wtilde u_{1}(\eta)-\wtilde u_{2}(\eta)))
\end{align*}%
where $\wtilde u_{a}=\ha(\wtilde u_{1}+\wtilde u_{2})$.

Since $|Q_{kk}-I|<C\epsilon_{0}$, for $\epsilon_{0}$ small enough (depending only on $n,k$) we
then have
\begin{align*}%
  \tag*{(3)}&|u_{a}(x)-u_{a}(x_{0})-(x-x_{0})Du_{a}(x_{0})| \\ %
  &\hskip.8in \le C\bigl(|\wtilde u_{a}(\xi)|+|\wtilde u_{a}(\eta)|+ %
  |(\wtilde u_{1}(\xi)-\wtilde u_{2}(\xi)) -(\wtilde u_{1}(\eta)-\wtilde u_{2}(\eta))|\bigr).
\end{align*}%
Notice also that by taking differences in~(1) we have
$$%
(\xi-\eta,\wtilde u_{1}(\xi)-\wtilde u_{2}(\eta))=(0,u_{1}(x)-u_{2}(x)){\cal{}Q}
$$%
whence
$$%
|\xi-\eta|<C|v(x)|\le C\epsilon_{0}\rho^{3/2} %
\leqno{(4)}
$$%
by~7.1.

Taking $\epsilon_{0}=\epsilon_{0}(n,k)$ small enough to ensure $C\epsilon_{0}<{1\over{}4}$, we
consider the following 2 cases with $\sigma=\rho^{3/2}$: Case~1 $B_{2\sigma}(\xi)\cap
{\cal{}K}_{\wtilde u}\neq \emptyset$ and Case~2 $B_{2\sigma}(\xi)\cap {\cal{}K}_{\wtilde u}=
\emptyset$.  In Case~1 we have by~7.1 (applied with $\wtilde u$ in place of $u$) and~(4) that
$|\wtilde u_{1}(\xi)-\wtilde u_{2}(\xi)|+|\wtilde u_{1}(\eta)-\wtilde u_{2}(\eta)|\le
C\epsilon_{0}\sigma^{3/2}\le C\epsilon_{0}\rho^{9/4}\le C\epsilon_{0}\rho^{2}$. In Case~2, for
$\epsilon_{0}=\epsilon_{0}(n,k)$ small enough, we also have $B_{3\sigma/2}(x)\cap
{\cal{}K}_{u}=\emptyset$ and hence there is a unique pair $u_{1},u_{2}$ of smooth single valued
functions on the ball $B_{3\sigma/2}(x)$ with $u|B_{3\sigma/2}(x)=\{u_{1},u_{2}\}$ and
corresponding smooth single valued $\wtilde u_{1},\wtilde u_{2}$ on $B_{\sigma}(\xi)$ so
that~(1),(3) hold. Then by 1-variable calculus along the line segment joining $\xi$ to $\eta$ we get
$$%
|(\wtilde u_{1}(\xi)-\wtilde u_{2}(\xi)) -(\wtilde u_{1}(\eta)-\wtilde u_{2}(\eta))| \le %
\sup_{B_{\sigma}(\xi)}|D(\wtilde u_{1}-\wtilde u_{2})|\sigma %
\leqno{(5)}
$$%
and by~7.1 (applied with $\wtilde u$ in place of $u$) we have $\sup_{B_{\sigma}(\xi)}|D(\wtilde
u_{1}-\wtilde u_{2})|\le C\epsilon_{0}\rho^{1/2}$, and so the right side of~(5) is $\le
C\epsilon_{0}\rho^{1/2}\sigma= C\epsilon_{0}\rho^{2}$.  Thus in both Case~1 and Case~2 we
conclude
$$%
|(\wtilde u_{1}(\xi)-\wtilde u_{2}(\xi)) -(\wtilde u_{1}(\eta)-\wtilde u_{2}(\eta))| \le %
C\epsilon_{0}\rho^{2},
$$%
and hence using this together with~7.2 (applied to $\wtilde u_{a}$ instead of $u_{a}$) on the right
of~(3) we obtain
$$%
|(u_{a}(x)-u_{a}(x_{0}))-(x-x_{0})Du_{a}(x_{0})|\le C\epsilon_{0}\rho^{2}
$$%
as claimed

\medskip

Finally we show that $u_{a}\in C^{1,1}$:

\begin{state}{\bf{}7.4 Theorem.} %
  $u_{a}\in C^{1,1}(B_{1/8})$ with $\sup_{B_{1/8}}|D^{2}u_{a}|\le C\epsilon_{0}$.
\end{state}%

{\bf{}Proof:} The proof is based on the estimates of Theorem~7.1, Lemma~7.3 and elliptic estimates
for single valued solutions. First note that by equation~5.11' and Theorem~7.1 we have
$$%
\ssum_{ij} a^{ij}_{\kappa\lambda}D_{i}D_{j}u_{a}^{\lambda} = F_{\kappa}\text{ with } %
|a^{ij}_{\kappa\lambda}-\delta_{ij}\delta_{\kappa\lambda}|\le C\epsilon_{0}, %
\text{ and }|F_{\kappa}|\le C\epsilon_{0}, \quad C=C(n,k),%
\leqno{\hbox{(1)}}
$$%
$\kappa=1,\ldots,k$, where we used Theorem~7.1 in checking that $|F_{\kappa}|\le C\epsilon_{0}$. 
It is then standard (the $L^{2}$ elliptic theory applied to the function
$u_{a}-u_{a}(x_{0})-(x-x_{0})\cdot Du_{a}(x_{0})$, which satisfies the same equation~(1)) that,
provided $\epsilon_{0}=\epsilon_{0}(n,k)$ is small enough,
$$%
\sint_{B_{\rho/2}(x_{0})} |D^{2}u_{a}|^{2} \le %
C\rho^{-4}\sint_{B_{\rho}(x_{0})}|u_{a}-u_{a}(x_{0})-(x-x_{0})\cdot Du_{a}(x_{0})|^{2} + %
C\epsilon_{0}^{2}\rho^{n} %
$$%
for $x_{0}\in B_{1/4}\cap {\cal{}K}_{u}$ and $\rho\in (0,1/4]$, and by virtue of Lemma~7.3 this
implies
$$%
\rho^{-n}\sint_{B_{\rho/2}(x_{0})} |D^{2}u_{a}|^{2} \le C\epsilon_{0}^{2}, %
\quad x_{0}\in B_{1/4}\cap {\cal{}K}_{u}, \,\rho\in (0,1/4], %
\leqno{\hbox{(2)}}
$$%
with $C$ independent of $x_{0}$ and $\rho$, and hence $|D^{2}u_{a}(x_{0})|\le C\epsilon_{0}$ for
a.e.\ $x_{0}\in {\cal{}K}_{u}$.  (Of course this latter statement is vacuous if ${\cal{}K}_{u}$ has
measure zero, which we ultimately show must be the case.)

Now suppose that $x\in B_{1/8}\setminus{\cal{}K}_{u}$, choose $x_{0}\in{\cal{}K}_{u}\cap B_{1/4}$
with $|x-x_{0}|=d(x)$, and let $\rho=|x-x_{0}|$. Then $u|B_{\rho}(x)$ can be represented as an
ordered pair of $C^{\infty}$ solutions of the minimal surface equation, each with gradient of length
$\le C\epsilon_{0}$, and 7.1 plus quasilinear elliptic estimates implies that the equation~5.10' can
be written
$$%
\Delta v^{\kappa}+b^{ij}_{\kappa\lambda}D_{i}D_{j}v^{\lambda}=0 %
\text{ on }B_{\rho}(x), %
\leqno{\hbox{(3)}}
$$%
with
$$%
\sup_{B_{\rho/2}(x)} \rho|D_{\ell}b^{ij}_{\kappa\lambda}|+ %
\rho^{3/2}[D_{\ell}b^{ij}_{\kappa\lambda}]_{1/2,B_{\rho/2}(x)}\le C\epsilon_{0}, %
\leqno{\hbox{(4)}}
$$%
By~(4), Schauder theory (for single-valued solutions) can be applied in~(3), giving
$$%
[D^{2}v]_{1/2,B_{\rho/4}(x)}\le C\rho^{-3/2}\epsilon_{0}{\sup}_{B_{\rho/2}(x)}|Dv|\le %
C\epsilon_{0} \rho^{-1} \text{ by~7.1.} %
$$%
Via another application of~7.1, this shows that the function $F_{\kappa}$ on the right of~(1)
satisfies $[F_{\kappa}]_{1/2,B_{\rho/2}(x)}\le C\epsilon_{0}\rho^{-1/2}$. Thus using~(1) in
combination with Schauder theory gives
$$%
\rho^{1/2}[D^{2}u_{a}]_{1/2,B_{\rho/4}(x)}\le C\epsilon_{0}. %
\leqno{\hbox{(5)}}
$$%
On the other hand by~(2) with $4\rho$ in place of $\rho$ we know there is a set of positive measure
in $B_{\rho/4}(x)$ with $|D^{2}u_{a}|\le C\epsilon_{0}$ and then~(5) gives
$\sup_{B_{\rho/4}(x)}|D^{2}u_{a}|\le C\epsilon_{0}$ and hence in particular $|D^{2}u_{a}(x)|\le
C\epsilon_{0}$ as required.

\section*{8\quad A frequency function for \boldmath{$v$} and the
dimension of \boldmath{${\mathcal K}_{u}$}.}
           
Using the key regularity results $u_{a} \in C^{1, 1}(B_{1})$ and $v \in C^{1, 1/2}(B_{1})$, we
can now establish the monotonicity of a frequency function for $v$, and use it to bound the size of
${\mathcal K}_{u}$.

We first want to show that the above regularity results make it possible to obtain a suitable
frequency function by directly modifying the work of Garofalo and Lin~\cite{GarL86} to handle the
present 2-valued setting and higher codimension. (Codimension $k>1$ implies that stationarity of
$G=\graph u$ puts us in an elliptic system setting rather than the single equation setting discussed
in~\cite{GarL86}.)

In view of the estimates~7.1, 7.4 we see that the equation 5.11 for $v^{\kappa}$ can be written in
the form
$$%
\ssum_{i,j=1}^{n}D_{i}\bigl(A^{ij}D_{j}v_{\kappa}\bigr)+ %
\ssum_{\ell=1}^{n}\ssum_{\lambda=1}^{k}E_{\kappa\lambda}^{\ell}D_{\ell}v^{\lambda}=0, %
\quad \kappa=1,\ldots,k, %
\leqno{\hbox{\bf 8.1}}
$$%
with $A^{ij}$ Lipschitz (real-valued) with small Lipschitz constant and $E_{\kappa}^{\ell\lambda}$
bounded:
$$%
|A^{ij}(x)-A^{ij}(y)|\le C\epsilon_{0}|x-y|, \quad |E_{\kappa\lambda}^{\ell}(x)| \le C \,\,\forall x,y\in
B_{1/2}, \leqno{\hbox{\bf 8.2}}
$$%
and $A^{ij},E_{\kappa}^{\ell\lambda}$ single-valued.  It is of crucial importance that we can thus
write the equation for $v$ as a system which is only weakly coupled (i.e.\ the top order part
$D_{i}(A^{ij}D_{j}v^{\kappa})$ involves application of the same scalar second order operator
$D_{i}(A^{ij}D_{j})$ to each component $v^{\kappa}$ of $v$).

We now claim that, in view of 8.1, 8.2 (which depend of course on the main regularity results~7.1
and~7.4), we can make a straightforward modification of the work of Garofalo \& Lin to establish a
key monotonicity result for the function $v$.

Before stating the result we need to set up some notation:

In view of~8.2 we can assume without loss of generality that
$$%
\det (A^{ij}) \equiv 1,
$$%
because otherwise we can replace $A^{pq}$ by $\det(A^{ij})^{-1/n}A^{pq}$ without changing the
form of the equation~8.1 and at the same time ensuring~8.2 still holds with $C=C(n,k,\alpha)$. 
Then as in~\cite{GarL86} we can view the operator $D_{i}(A^{ij}D_{j})$ in~8.1 as 
 the Laplacian with respect to the metric $\sum_{i,j}A_{ij}dx_{i}dx_{j}$ where 
$$%
(A_{ij})=(A^{ij})^{-1}.
$$%
The fact that $A_{ij}$ is Lipschitz is not sufficient to introduce normal coordinates, but
following~\cite{AKS61} we can first multiply by the Lipschitz conformal factor
$$%
\eta(x)=A^{\ell m}(x)(x_{\ell}/r) (x_{m}/r), \quad
r=(\ssum_{j}x_{j}^{2})^{1/2},
$$%
to give a new metric
$$%
\wtilde A_{ij} =\eta A_{ij} = A^{\ell m}(x)(x_{\ell}/r) (x_{m}/r) A_{ij}.
$$%
Notice that indeed $\eta$ is Lipschitz, because $A^{ij}=\delta_{ij}+a_{ij}$ with $a_{ij}$ Lipschitz and
$a_{ij}(0)=0$, and hence $\eta=1+\sum_{i,j}a_{ij}(x_{i}/r)(x_{j}/r)$ which implies $|D_{\ell}\eta|\le
C\sum_{i,j}(|D_{\ell}a_{ij}|+r^{-1}|a_{ij}|)\le C$. Thus the equation~8.1 can be written
$$%
\ssum_{i,j=1}^{n}D_{i}\Bl(\wtilde A^{ij}D_{j}v_{\alpha}\Br)+ %
\ssum_{\ell=1}^{n}\ssum_{\beta=1}^{k}\wtilde E_{\alpha}^{\ell\beta}D_{\ell}v_{\beta}=0, \quad
\alpha=1,\ldots,k,
$$%
where
$$%
\wtilde A^{ij} = (A^{\ell m}(x)(x_{\ell}/r) (x_{m}/r))^{\fr n-2/2}A^{ij}.
$$%
Thus working with $\wtilde g_{ij}$ instead of $g_{ij}$ involves merely multiplying the principal
coefficients of the equation by the Lipschitz function $\eta^{\fr n-2/2}$, where
$\eta=\tsum_{i,j}A^{ij}(x)(x_{i}/r) (x_{j}/r)$, and this does not change the form of the equation (or
the boundedness of the coefficients of the first order terms). But it has the advantage (as proved
in~\cite{AKS61}) that, provided we have~8.2 with $\epsilon_{0}=\epsilon_{0}(n,k)$ small enough,
we can make a bilipschitz change of coordinates $y=\Gamma(x)$ with $\Gamma(0)=0$,
$$%
\begin{aligned}%
  \Gamma(B_{1})&\supset B_{1/2} ,\quad \rho_{0}=\rho_{0}(n), \\ %
  \wtilde A_{ij}dx_{i}dx_{j}&=\what g_{ij}(y) dy_{i} dy_{j} %
\end{aligned}%
$$%
such that the new coordinates have some of the key properties enjoyed by normal coordinates,
including
$$%
\ssum_{j}\what g_{ij}y_{j} = y_{i} \text{ on }B_{1/2}, \quad i=1,\ldots,n,
$$%
and also having the property that the radial derivatives $D_{r}\what g_{ij}\equiv
\ssum_{j}|y|^{-1}y_{j}D_{y_{j}}\what g_{ij}$ are bounded:
$$%
|D_{r}\what g_{ij}| \le C, \quad C=C(n,k), i,j=1,\ldots,n.
$$%
Then, writing $\what A^{ij}=\sqrt{\what g}\hskip1.5pt\what g^{ij}$, the equation~8.1 transforms to a
new equation of the same form, Viz.,
$$%
\ssum_{i,j=1}^{n}D_{i}\Bl(\what A^{ij}D_{j}\what v^{\kappa}\Br)+ %
\ssum_{\ell=1}^{n}\ssum_{\lambda=1}^{k}\what E_{\kappa\lambda}^{\ell}%
D_{\ell}\what v^{\lambda}=0, \quad \kappa=1,\ldots,k, %
\leqno{\hbox{\bf 8.3}}
$$%
where $\what v=v\circ \Gamma^{-1}$,
$$%
\what A^{ij}(0)=\delta_{ij},\quad |D_{r}\what A^{ij}|\le C\,, %
\footnote[3]{\,We emphasize that the construction of~\cite{AKS61} only ensures boundedness of the
radial derivative $D_{r}\what A^{ij}$, not the tangential derivatives.}%
\quad |\what E_{\alpha}^{\ell\beta}| \le C,%
\leqno{\hbox{\bf 8.4}}
$$%
and where we now also have
$$%
\begin{aligned}%
  &\ssum_{j}\what A^{ij}y_{j} \equiv \mu y_{i},\,\,y\in B_{1/2}, %
  \quad \mu=\sqrt{\what g}\quad i=1,\ldots,n,\\ %
  &|D_{r}\mu| \le C, \quad C^{-1}\le \mu \le C, \quad C=C(n,\alpha). %
\end{aligned}%
\leqno{\hbox{\bf 8.5}}
$$%

\begin{state}{\bf{}8.6 Lemma.} %
  There is $\epsilon_{0}=\epsilon_{0}(n,k,\alpha)$ such that if $0\in{\cal{}K}_{u}$, if~1.2, 1.6 hold,
and if $\what v=v\circ \Gamma^{-1}$,  $\what A^{ij}$, and $\mu$
are as above, then the modified frequency function $\what N_{\what v}(\rho)$ defined by
$$%
\what N_{\what v}(\rho) = %
\Fr \rho^{2-n}\mint_{\partial B_{\rho}}\mu \what v\cdot\what v_{r} \, d{\cal{}H}^{n-1}
/{\rho^{1-n}\mint_{\partial B_{\rho}} \mu\, |\what v(y)|^{2}\,d{\cal{}H}^{n-1}(y)}
$$%
has the property that $(\exp C\rho)\what N_{\what v}(\rho)$ is increasing as a function of $\rho\in
(0,\rho_{0}]$, where $C=C(n,k,\alpha),\,\rho_{0}=\rho_{0}(n,k,\alpha)\in (0,\fr 1/2]$, and
furthermore we can choose
$\theta=\theta(n,k,\alpha)\in (0,1)$ so that 
 for each $\beta>\what N_{\what v}(\rho_{0})$
we have the fixed lower bounds {\belowdisplayskip-1pt
$$%
\sigma^{1-n}\int_{\partial
B_{\sigma}}|\what v|^{2}\,d{\cal{}H}^{n-1} \ge  (\sigma/\rho)^{2\beta}\rho^{1-n}\int_{\partial
B_{\rho}}|\what v|^{2}\,d{\cal{}H}^{n-1},\quad 0<\sigma\le \rho\le \min\{\rho_{0},\theta(\fr
\beta/{\what N_{\what v}(\rho_{0})}-1)\}. 
$$}%
\end{state}

{\bf{}8.7 Remark:} Notice that the modified frequency $\what N_{\what v}$ is of similar order to
the usual frequency function $N_{\what v}(\rho)={\cal{}D}(\rho)/{\cal{}H}(\rho)$, where
${\cal{}D}(\rho)=\rho^{2-n}\int_{B_{\rho}}\what A^{ij}\what v_{i}\cdot\what v_{j}$,
${\cal{}H}(\rho)=\rho^{1-n}\int_{\partial B_{\rho}}\mu |\what v|^{2}\,d{\cal{}H}^{n-1}$, by virtue of
the fact that
$$%
(1-C\rho){\cal{}D}(\rho) \le {\cal{}I}(\rho) \le (1+C\rho){\cal{}D}(\rho)
$$%
for sufficiently small $\rho$, where ${\cal{}I}(\rho)$ is the quantity on the top line of $\what
N_{\what v}(\rho)$ in~8.6 above, because~8.3 implies
$$%
{\cal{}D}(\rho) = {\cal{}I}(\rho)+ \rho^{2-n}\int_{B_{\rho}}R(\what v)\cdot\what v
$$%
with $R(\what v)$ an $\R^{k}$-valued function such that $|R(\what v)|\le C(|\what v|+|D\what
v|)$, and~6.5 is applicable with $w=\what v$. (Only the fact that $w$ satisfies an equation of the
form of~8.3 and the fact that $w$ was Lipschitz, rather than $C^{1,\alpha}$, was used to prove~6.5.)

{\bf{}Proof of~8.6.}  First note that (by applying the operator $\tsum y_{j}D_{y_{j}}$ to equation 8.3
and using~8.4) we see
$$%
\tsum_{j}y_{j}D_{y_{j}}\what v_{\lambda}\in W^{1,2}(B_{1/2}\setminus {\cal{}K}_{\what v}), \quad
\lambda=1,\ldots,k.
 $$%
 Also by multiplying by $\gamma_{\delta}(\what v)$ in~8.3, with $\gamma_{\delta}(t)=\sgn t\,
\max\{|t|-\delta,0\}$, and integrating by parts we see that (after letting $\delta\downarrow0$)
$$%
\sint_{B_{\rho}} |D\what v|^{2} \le C\rho\sint_{B_{\rho}}(|\what v|^{2}+ |D\what v|^{2})
$$%
if $\what v\equiv\{0,0\}$ on $\partial B_{\rho}$ and hence for small enough $\rho$ we use~6.5
(with $w=\what v$) to conclude that $\what v\equiv 0$ on $B_{\rho}$. Thus the modified frequency
function $\what N_{\what v}(\rho)$ is well-defined for $\rho\in (0,\rho_{0}]$ for suitable
$\rho_{0}=\rho_{0}(n,k,\alpha)$.

 In view of~8.3--8.5 we can now check the two basic identities (analogous to the identities in
Remark~2.3(2) in the special case $\zeta_{j}(x)\equiv x_{j}$) :
$$%
\mathcal{D} \equiv  \rho^{2-n}\int_{B_{\rho}} \what A_{kl} \what v_k\cdot \what v_l =   %
\rho^{2-n}\int_{\partial B_{\rho}} \mu  \what v\cdot \what v_r + 
                                    \rho^{2-n} \int_{B_{\rho}} R(\what v) \cdot\what v, 
\leqno{\hbox{\bf 8.8}}
$$%
\begin{align*}%
  & (n-2)\rho^{1-n}\int_{B_{\rho}}\what A^{kl}\what v_k\cdot \what v_l -  %
       \rho^{2-n} \int_{\partial B_{\rho}} \what A^{kl}\what v_k\cdot \what v_l =  %
                -\rho^{2-n} \int_{\partial B_{\rho}}2 \mu  |\what v_r|^2 \\ %
& \hskip2.5in - \rho^{1-n} \int_{B_{\rho}}
r(\what A^{kl,\,r}\what v_k\cdot\what v_l - 2 R(\what v)\cdot\what v_r),   %
\end{align*}%
i.e.\
\begin{align*}%
\tag*{\hbox{\bf 8.9}} &\mathcal{D}' = 
      \frac{d}{d\rho} \Big( \rho^{2-n} \int_{B_{\rho}} \what A^{kl} \what v_k\cdot \what v_l\Big) \\ %
 & \;\;\; = \rho^{2-n} \int_{\partial B_{\rho}}  %
                                 2\mu |\what v_r|^2 + \rho^{1-n} \int_{B_{\rho}} %
r(\what A^{kl}_{r}\what v_k\cdot \what v_l - 2 R(\what v)\cdot\what v_r) %
\end{align*}%
where subscripts denote partial derivatives, $|R(\what v)|\le C(|\what v|+|D\what v|)$ and
$\mu=\sqrt{\what g}$ as above. Using these identities we can now establish the required
monotonicity as in~\cite[pp.358-364]{GarL87} with $rf(r)\equiv $const.  (Notice we use~8.7, 8.8
and 8.9 in lieu of the corresponding inequalities/identities on p.358 of~\cite{GarL87} and we do not
need to refer back to~\cite{GarL86} as is done in~\cite{GarL87}).

To prove the last part of~8.6 observe that by the monotonicity of $\what N_{\what v}(\rho)$ and
the inequalities of Remark~8.7 we have $N_{\what v}(\rho)\le (1+C\rho) N_{\what v}(\rho_{0})$ for
$\rho\in (0,\rho_{0}]$. Also, by~8.8 and the fact that~6.5 is applicable with $w=\what v$ as
discussed in~8.7, we have ${\cal{}H}'(\rho)\le 2(1+C\rho) {\cal{}D}(\rho)$, and then we have
$\rho{\cal{}H}^{'}(\rho)/{\cal{}H}(\rho)\le 2(1+C\rho)\what N_{\what v}(\rho_{0})\le 2\beta$ if
$\rho\le\rho_{1}$, where $\rho_{1}$ is a suitably small multiple (depending in $n,k,\alpha)$ of
$\fr \beta/{\what N_{\what v}(\rho_{0})}-1$.
Hence we can integrate (as in Remark~2.3(4)) to conclude the stated bounds in the last part of~8.6.

\begin{state} {\bf 8.10 Theorem.} %
  Suppose $v \not\equiv 0$. Then the Hausdorff dimension of\/ ${\mathcal K}_{v}$ is at most
$(n-2)$. Furthermore, either ${\mathcal B}_{u} = \emptyset$ or the Hausdorff dimension of
$\,{\mathcal B}_{u}$ is equal to $(n-2)$ and the $(n-2)$-dimensional Hausdorff measure of
$\,{\mathcal B}_{u}$ is positive.
\end{state}

{\bf Proof:} For the first claim of the theorem, it suffices to show that ${\mathcal H}^{s} \,
({\mathcal K}_{v} \cap B_{1/2}) = 0$ for every $s > n-2$. So fix $s > n-2$  and suppose that
${\mathcal H}^{s} \, ({\mathcal K}_{v} \cap B_{1/2}) > 0$. Let $\mu_{s}$ be the outer measure on
$\R^{n}$ as defined in~\S4 and, also as in~\S4, let $z \in {\mathcal K}_{v} \cap B_{1/2}$ be a point
of positive upper density with respect to $\mu_{s}$.  Thus there exists a sequence of positive
numbers $\sigma_{j} \to 0$ such that
$$%
\lim_{j\to\infty}\,\sigma_{j}^{-s}\mu_{s}\,({\mathcal K}_{v}\cap B_{\sigma_{j}}(z)) >0. %
\leqno{(1)}
$$%
Let $v_{z}(x)=v(z+x)$. By~7.1 and~7.4 we have the interior $C^{1, 1/2}$ and $W^{2, 2}$
estimates
$$%
\sup_{B_{\theta\sigma}(y)}\, |v_{z}| + \sigma\sup_{B_{\theta\sigma}(y)} |Dv_{z}| + %
\sigma^{3/2}\hskip-.3in\sup_{x_{1}, x_{2} \in B_{\theta\sigma}(y), x_{1} \neq x_{2}} \, %
\frac{|Dv_{z}(x_{1}) - Dv_{z}(x_{2})|}{|x_{1} - x_{2}|^{1/2}} \leq %
C\left(\sigma^{-n}\int_{B_{\sigma}(y)}|v_{z}|^{2}\right)^{1/2} %
\leqno{(2)}
$$%
and
$$%
\sigma^{4-n}\int_{B_{\theta\sigma}(y)}|D^{2}v_{z}|^{2} \leq
C\sigma^{-n}\int_{B_{\sigma}(y)}|v_{z}|^{2} \leqno{(3)}
$$%
for each $\sigma\in (0,1/4)$, $y\in B_{1/4}(z)$ and each $\theta \in (0, 1)$, where $C = C(n, k,
\theta,\alpha)$. It follows directly from these estimates that if we let $v_{z, \sigma_{j}}(x) =
\sigma_{j}^{-n/2}\|v_{z}\|_{L^{2}(B_{\sigma_{j}})}^{-1} v_{z}(\sigma_{j}x)$ for $x \in B_{1}$, then
after passing to a subsequence, we have that $v_{z, \sigma_{j}} \to \varphi$ for some 2-valued,
symmetric function $\varphi \in C^{1, 1/2}(\R^{n}) \cap W^{2, 2}_{\text{loc}} (\R^{n})$, where the
convergence is in $C^{1,\beta}(B_{\rho})$ for every $\rho >0$ and $\beta<\ha$, guaranteeing also
that $\varphi$ is harmonic. Also, as in the proof of~4.1, we have $\mu_{s}({\cal{}K}_{\varphi}\cap
B_{1})>0$.

We claim that $\varphi$ is not identically zero. Indeed, by taking $\sigma=\theta\rho$ in the
inequality in the last part of~8.6 and integrating with respect to $\rho$, we have, for each $\theta\in
(0,1)$, $\int_{B_{\theta\rho}}|\what v|^{2}\ge C(n,k,\alpha,\theta)\int_{B_{\rho}}|\what v|^{2}$. 
Since the change of coordinates $\Gamma$ is bilipschitz this gives $\int_{B_{\rho/2}}|v|^{2}\ge
C(n,k,\alpha)\int_{B_{\rho}}|v|^{2}$, whence $\int_{B_{1/2}}|v_{z,\sigma_{j}}|^{2}\ge C$ for some
fixed $C>0$ (independent of $j$) and $\varphi$ is non-zero as claimed.  

Thus $\varphi$ is $C^{1,\alpha}$ harmonic, not identically zero, and
${\cal{}H}^{s}({\cal{}K}_{\varphi})>0$ for some $s>n-2$, which contradicts Lemma~4.1.

For the remaining assertion, note that if ${\mathcal H}^{n-2} \, ({\mathcal B}_{u}) = 0$ then by the
result of the Appendix $B_{1} \setminus {\mathcal B}_{u}$ is simply connected, so that $\left. 
u\right|_{B_{1} \setminus {\mathcal B_{u}}} = \{u_{1}, u_{2}\}$ for a pair of smooth functions
$u_{1}, u_{2} \, : \, B_{1} \setminus {\mathcal B}_{u} \to \R^{k}$ each solving the minimal surface
system. Since ${\mathcal H}^{n-1} ({\mathcal B}_{u}) = 0$ and $u_{1},u_{2}$ are
$C^{1,\alpha}(B_{1})$ it follows that in fact $u_{1},u_{2}$ are weak $C^{1,\alpha}(B_{1})$ (hence
smooth strong) solutions of the minimal surface system, which implies that ${\mathcal B}_{u} =
\emptyset$.

\bigskip

Combining the codimension~1 case of Theorem~8.10 with the main regularity theorem
of~\cite{Wic08}, we obtain the following:

\begin{state}{\bf{}8.11 Theorem.} %
  Let $V$ be an $n$-dimensional stationary integral varifold in an open set $U\subset \R^{n+1}$
arising as the weak limit of a sequence of stable minimal hypersurfaces $M_{j}$ of $U$ with
$M_{j}$ immersed away from a closed set $K_{j}$ of locally finite $(n-2)$-dimensional Hausdorff
measure for each $j$. Then the set of points $z \in \spt \|V\|$ where $V$ has a multiplicity 2
tangent plane but $\spt \|V\|$ is not a smooth embedded submanifold near $z$ has Hausdorff
dimension at most $(n-2).$ In particular, the set of multiplicity 2 branch points of $\,V$ (i.e. the set
of points $z \in \spt \|V\|$ with the property that $V$ has a multiplicity 2 tangent plane at $z$, but
for no $\sigma>0$ is $\spt \|V\| \cap B_{\sigma}^{n+1}(z)$ equal to a smooth immersed hypersurface
of $B_{\sigma}^{n+1}(z)$) is either empty or has Hausdorff dimension equal to $(n-2)$ and locally
positive $(n-2)$-dimensional Hausdorff measure.
\end{state}%

\newpage

\section*{Appendix: A simple connectivity lemma.}

The following lemma is presumably well known, but we include it here for the convenience of the
reader since we have not found it in the literature.

\begin{state}{\bf{}Lemma.} %
  If\/ $\Gamma$ is a closed subset of\/ $\R^{n}$ with ${\cal{}H}^{n-2}(\Gamma)=0$, then
$B_{1}\setminus \Gamma$ is simply connected.
\end{state}

{\bf{}Proof:} We use induction on $n\ge 2$. In case $n=2$ the result is trivial because
$\Gamma=\emptyset$ in this case. So assume that $n\ge 3$ and that the result holds with $n-1$ in
place of $n$. $\ovl{B}_{1}$ is bilipschitz homeomorphic to $\ovl{Q}_{1}$, where $Q_{1}$ is the
cylinder $B^{n-1}_{1}\times (-1,1)$, and this homeomorphism takes the compact set
$\Gamma\cap\ovl{B}_{1}$ to the compact set $\widetilde\Gamma$, where
${\cal{}H}^{n-2}(\widetilde\Gamma)=0$, so it suffices to prove $Q_{1}\setminus
\widetilde\Gamma$ is simply connected.

By the ``rough coarea inequality'' (\cite[2.10.25]{Fed69}) we can pick $y_{0}\in (-1,1)$ such that
$$%
{\cal{}H}^{n-3}((\R^{n-1}\times\{y_{0}\}) \cap \widetilde\Gamma)=0. %
\leqno{(1)}
$$%
Let $\gamma:S^{1}\to\R^{n}$ be an arbitrary smooth closed curve contained in
$Q_{1}\setminus\widetilde\Gamma$, let
$$%
\delta=\min\{|x-y|:x\in \gamma,\,y\in\widetilde\Gamma\},
$$%
and let $P$ denote the projection $(x^{1},\ldots,x^{n-1},x^{n})\mapsto (x^{1},\ldots,x^{n-1},0)$ of
$\R^{n}$ onto $\R^{n-1}\times\{0\}$. Then $P(\widetilde\Gamma)\times\R$ is a compact subset of
$\R^{n}$ with $(n-1)$-dimensional Hausdorff measure zero and hence its orthogonal projection
onto any given hyperplane in $\R^{n}$ also has $(n-1)$-dimensional Hausdorff measure zero.  Thus
if $\eta\in S^{n-1}$ and $L$ is the hyperplane orthogonal to $\eta$, and if $P_{\eta}$ is orthogonal
projection of $\R^{n}$ onto $L$, then
$A_{\eta}=P_{\eta}^{-1}(P_{\eta}(P(\widetilde\Gamma)\times\R))$ is a set of
${\cal{}H}^{n}$-measure zero and every line $\ell_{y}(\eta)=\{y+\tau\eta:\tau\in\R\}$ with $y\notin
A_{\eta}$ is disjoint from $A_{\eta}\cup(P(\widetilde\Gamma)\times\R)$.  It is then elementary to
construct an approximation $\widetilde\gamma:S^{1}\to\R^{n}$ to $\gamma$ with $\max_{t\in
S^{1}}|\widetilde\gamma(t)-\gamma(t)|<\delta$ and $\widetilde\gamma\subset
Q_{1}\setminus(P(\wtilde\Gamma)\times\R)$.  (We can construct such a $\widetilde\gamma$ by
first taking a polygonal approximation to $\gamma$ with edge segments $s_{1},\ldots,s_{N}$ parallel
to $\eta_{1},\ldots,\eta_{N}\in S^{n-1}$ respectively, and then composing with a translation $\tau$
of $\R^{n}$ chosen to ensure that none of the translated segments $\tau(s_{j})$ have endpoints in
$A_{\eta_{j}}$, $j=1,\ldots,N$.)

Then $\gamma$ can be homotopied in $Q_{1}\setminus\widetilde\Gamma$ to
$\widetilde\gamma$ via the homotopy $\gamma_{s}(t)=s\widetilde\gamma(t)+(1-s)\gamma(t)$,
$s\in [0,1]$), and $\wtilde\gamma$ can be homotopied in
$Q_{1}\setminus(P(\widetilde\Gamma)\times\R)$ to $\widetilde\gamma_{1}\subset
Q_{1}\cap(\R^{n-1}\times\{y_{0}\})$ via the homotopy $\widetilde\gamma_{s}(t)=
(\widetilde\gamma_{1}(t),\ldots,\widetilde\gamma_{n-1}(t),sy_{0}+(1-s)\widetilde\gamma_{n}(t)),
\,s\in [0,1]$.  Finally by~(1) and the inductive hypothesis $\widetilde \gamma_{1}$ can be
homotopied in $Q_{1}\cap(\R^{n-1}\times\{y_{0}\})\setminus\widetilde\Gamma$ to a point, and
the proof is complete.

\providecommand{\bysame}{\leavevmode\hbox to3em{\hrulefill}\thinspace}
\providecommand{\href}[2]{#2}

\bigskip

$$
\hskip-.2in\vbox{\hsize2.5in\obeylines\parskip-4pt 
  \small 
Leon Simon 
Mathematics Department 
Stanford University 
Stanford CA 94305, USA
lms@math.stanford.edu} 
\vbox{\hsize3in\obeylines \parskip-4pt \small Neshan Wickramasekera
DPMMS 
University of Cambridge 
Cambridge CB3 0WB, United Kingdom
N.Wickramasekera@dpmms.cam.ac.uk}
$$

\end{document}